\documentclass[10pt]{article}
\usepackage{amsmath}
\usepackage{amssymb}
\usepackage{amsthm}
\usepackage{graphicx}
\theoremstyle{plain}
\newtheorem{thm}{Theorem}[section]
\newtheorem{cor}[thm]{Corollary}
\newtheorem{lem}[thm]{Lemma}
\theoremstyle{definition}
\newtheorem{defn}{Definition}[section]
\newtheorem{claim}{Claim}

\theoremstyle{remark}
\newtheorem{rem}{Remark}[section]
\begin{document}
\newlength{\figwidth}
\setlength{\figwidth}{\textwidth}

\newlength{\fighalfwidth}
\setlength{\fighalfwidth}{0.48\textwidth}

\newlength{\fighalfheight}
\setlength{\fighalfheight}{0.2\textheight}

\newlength{\figthirdwidth}
\setlength{\figthirdwidth}{0.3\textwidth}
\title{Pitchfork bifurcations of invariant manifolds}

\author{Jyoti Champanerkar and Denis Blackmore\\
{\small(Department of Mathematical Sciences, New Jersey Institute of Technology)}\\
email: jac4@njit.edu}

\date{}

\maketitle
\begin{abstract}
A pitchfork bifurcation of an $(m-1)$-dimensional invariant submanifold of a
dynamical system in $\mathbb{R}^m$ is defined analogous to that in $\mathbb{R}$. Sufficient
conditions for such a bifurcation to occur are stated and existence of the bifurcated
manifolds is proved under the stated hypotheses. For discrete dynamical systems, the
existence of locally attracting manifolds $M_+$ and $M_-$, after the bifurcation has taken
place is proved by constructing a diffeomorphism of the unstable manifold $M$. For
continuous dynamical systems, the theorem is proved by transforming it to the discrete case.
Techniques used for proving the theorem involve differential topology and analysis.
The theorem is illustrated by means of a canonical example.
\end{abstract}

\section{Introduction}
Pitchfork bifurcations bear the name due to the fact that the bifurcation diagram for a
one-parameter family in $\mathbb{R}$ looks like a pitchfork. Pitchfork bifurcation for a
fixed point in $\mathbb{R}$ has been widely studied. In $\mathbb{R}$, sufficient conditions
for the occurrence of a pitchfork bifurcation of a non-hyperbolic fixed point are stated for
instance in \cite{Rasband, Wigg2}. A generalization of the result in $\mathbb{R}$ is given by
Sotomayor's theorem \cite{Perko} for a pitchfork bifurcation of a fixed point in $\mathbb{R}^n$.
Another generalization for pitchfork bifurcations is that for a periodic orbit \cite{Perko}.
Analytical discussions of pitchfork (or pitchfork type) bifurcations can be found for
particular classes of dynamical systems, e.g.,\cite{Glen}, where  a quasi-periodically forced map
is studied. Interesting numerical analyses of pitchfork bifurcations can be found
in \cite{OsWiGlFe, Stur}.

An algorithm to compute invariant manifolds of
equilibrium points and periodic orbits is presented in \cite{KrOs}. It is important to study
invariant manifolds in order to know the global dynamics of a system. The classical
pitchfork bifurcation concerns a fixed point (invariant codimension-1 submanifold) on
the real line. From a mathematical viewpoint, it is therefore natural and important to
investigate higher dimensional extensions of this theorem to invariant codimension-1
submanifolds of a Euclidean m-space. Accordingly we ask under what conditions
an invariant manifold of a (discrete or continuous) dynamical system undergoes
a pitchfork bifurcation? A fairly complete answer to this question is provided in
this work. We give sufficient conditions for the occurrence of a pitchfork bifurcation of a
compact, boundaryless, codimension-$1$, invariant manifold in $\mathbb{R}^m$. We obtain
readily verifiable criteria for identifying such bifurcations, and illustrate the use of
these criteria in an example. Techniques used for proving the theorem involve differential
topology and analysis and are adapted from Hartman \cite{Hartman}, Hirsch et al. \cite{HiPuSh}
and Shub \cite{Shub}.
\section{Definitions}
Let $M$ be a codimension-$1$, compact, connected, boundaryless submanifold of $\mathbb{R}^m$.
By the Jordan-Brouwer separation theorem \cite{GuPo}, $M$ divides $\mathbb{R}^m \backslash
M$ into an outer unbounded region and an inner bounded region. We shall study
$\mathcal{C}^1$ functions $F:U \times (-a,a) \rightarrow \mathbb{R}^m$, where $U$ is an open
neighborhood of $M$ in $\mathbb{R}^m$ and $(-a,a)$, $a>0$, is an open symmetric interval of
real numbers. It shall be assumed in the sequel that each of the maps $F_{\mu}:U\rightarrow
\mathbb{R}^m$, $|\mu|<a$, is a $\mathcal{C}^1$ diffeomorphism and that $M$ is
$F_{\mu}$-invariant, i.e. $F_{\mu}(M) =M$.
\begin{defn}
With $M$ and $F_{\mu}$ as above, we say that $F_{\mu}$ is side-preserving if
for every $x$ in the inner bounded region, $F_{\mu}(x)$ also lies in the inner region.
\end{defn}
\begin{defn}
With $M$ and $F_{\mu}$ as above, we say that $F_{\mu}$ is side-reversing if
for every $x$ in the inner bounded region, $F_{\mu}(x)$ lies in the outer unbounded region.
\end{defn}
Note that if $F_{\mu}$ is a diffeomorphism of a neighborhood of $M$ and leaves $M$ invariant, then
$F_{\mu}$ is either side-preserving or side-reversing. Observe also that in the case that
$F_{\mu}$ is side-reversing, it is not possible for $U=\mathbb{R}^m$.
Analogous to the definition
of a pitchfork bifurcation in $\mathbb{R}$, we define a
pitchfork bifurcation of invariant manifolds in $\mathbb{R}^m$ as follows.
\begin{defn}
Consider a discrete dynamical system in $\mathbb{R}^m$ given by $x_{n+1} = F_{\mu}(x_n)$.
Let $M$ be an invariant manifold for all $\mu \in (-a,a)$. If $0 \leq \mu_0  < a$ is such
that $M$ is locally attracting (repelling) for $\mu < 0$, $M$ is locally repelling
(attracting) for $\mu
> \mu_0$ and in addition two locally attracting (repelling) $F_{\mu}$-invariant diffeomorphic copies of
$M$, viz., $M_-$ and $M_+$ appear in a small neighborhood of $M$
for $\mu > \mu_0$, then we say that $M$ has undergone a pitchfork bifurcation at $\mu_0$.
\end{defn}
In the definition above, it does not matter what happens in the interval $(0, \mu_0)$. It is
typically assumed that the interval $(0, \mu_0)$ is small. In $\mathbb{R}$, $0$ coincides
with $\mu_0$, since the manifold under consideration is just a single point. But for higher
dimensions, not all points on the invariant manifold may undergo a change in stability at the same
value of the parameter $\mu$. When $\mu \geq \mu_0$, all the points have changed stability and two
new invariant, diffeomorphic copies of the original manifold (of opposite stability) appear.
\section{Pitchfork bifurcation theorem for discrete dynamical system}
\label{sec:PBdiscrete}
We consider (one-parameter families of) discrete dynamical systems given by
\begin{equation} \label{eq:discsystem}
x_{n+1}=F(x_n,\mu)
\end{equation}
where, $x_n \in \mathbb{R}^m$ for every $n \in \mathbb{N}$ and
$\mu \in (-a,a) \subseteq \mathbb{R}$. Additional properties of $F(\cdot ,\mu)$, also denoted as
$F_{\mu}(\cdot)$ are
described later. Let $M$ be a compact, connected, boundaryless, codimension-$1$, $\mathcal{C}^1$
submanifold of $\mathbb{R}^m$, which is $F_{\mu}$-invariant $\forall \mu \in (-a, a)$. Denote a tubular
neighborhood of $M$ as $N(\alpha)=\{x \in \mathbb{R}^m : d(x,M) \leq \alpha, \alpha >0 \}$, where
$d$ is the standard Euclidean distance function. Assume that
$\alpha$ is sufficiently small so that the $\epsilon$-neighborhood theorem \cite{GuPo}
can be applied and $N(\alpha) \subset U$. This means that every element $x \in N(\alpha)$
can be uniquely represented as $x=(r,y)$ where $y=\pi(x) \in M$ is the point on $M$ closest
to $x$ and $r \in [ -\alpha,\alpha ]$ is the signed distance in the outward normal direction
between $x$ and $M$. We also assume that $F_{\mu}(N(\alpha)) \subset N(\alpha)$. This
enables us to write $F_{\mu}$ in component form as $F_{\mu}=(f_{\mu},g_{\mu})$ where
$f_{\mu}:N(\alpha) \rightarrow \mathbb{R}$, $g_{\mu}=\pi \circ F_{\mu}: N(\alpha) \rightarrow
M$, and $f_{\mu}(x)$ is the signed distance from  $g_{\mu}(x)$ to $F_{\mu}(x)$.
Observe that $F_{\mu}^{-1}(F_{\mu}(N(\alpha))) = N(\alpha)$, so that $F_{\mu}^{-1}$ can be
written in component form as $F_{\mu}^{-1}=(\hat{f_{\mu}},\hat{g_{\mu}})$.

We shall use standard notation for derivatives and partial derivatives of functions. For
example, the derivative of $F_{\mu}:N(\alpha) \rightarrow \mathbb{R}^m$ will be denoted by
$DF_{\mu}$, and represented as the usual $m\times m$ Jacobian matrix
\begin{displaymath}
DF_{\mu}(r,y)=\left[
\begin{array}{cc} D_rf_{\mu} (r,y) &  D_yf_{\mu}(r,y)\\
D_rg_{\mu}(r,y)& D_yg_{\mu}(r,y)
\end{array}\right],
\end{displaymath}
where the entries are submatrices representing the partial derivatives such as
\begin{displaymath}
D_rf_{\mu}(r,y)=\frac{\partial f_{\mu}}{\partial r}(r,y) {\rm \quad and \quad}
D_yg_{\mu}(r,y)=\left[\frac{\partial g_{{\mu}_i}}{\partial y_j}\right]_{(m-1)\times (m-1)}.
\end{displaymath}
We use $|\cdot |$ for the Euclidean norm of an element of a Euclidean space or the
associated norm of a linear mapping (matrix) between Euclidean spaces. The symbol $\| \cdot
\|$ denotes the supremum norm of a function taking values in a Euclidean space or in a space
of linear transformations of Euclidean spaces taken over an appropriate set, which is
sometimes indicated as a subscript of the norm.

If the rate of change (with respect to $r$) in the normal component $f_{\mu}$ in the radial
direction $r$ is strictly less than $1$ in absolute value, the manifold $M$ will be locally
attracting. This is stated mathematically in statement $(ii)$ of Theorem \ref{thm:PBthethm}, which
follows.
Similarly to have $M$ locally repelling, we require that $|D_rf_{\mu}|$ be greater than one
as in statement $(iii)$. Statements $(iv)$ and $(v)$ describe locally attracting properties in a
neighborhood away from $M$, which is where our new bifurcated manifolds $M_-$ and $M_+$ will
reside. Properties $(vi)$ and $(vii)$ are obtained analytically and are needed in order to
establish the existence of manifolds $M_-$ and $M_+$ as graphs of a fixed point (Lipschitz
function) in a Banach space. The last hypothesis (statement $(viii)$) provides boundedness and
equicontinuity properties that enable us to bootstrap Lipschitz homeomorphisms of $M$
with $M_+$ and $M_-$ up to $\mathcal{C}^1$ diffeomorphisms.
These ideas shall become clear as the proof
unfolds and after the remarks following the proof.
\begin{thm}\label{thm:PBthethm}
With $F_{\mu}$ and $M$ as above, suppose that the following statements hold.
\begin{enumerate}
\item $F_{\mu}$ is side-preserving for every $\mu \in (-a,a)$.

\item $\underset{(r,y) \in N(\alpha)}{\sup}|D_rf_{\mu}(r,y)|=\|D_rf_{\mu}\|_{N(\alpha)} <
1$ for every $\mu \in (-a,0)$.

\item $\exists$ \quad $0< \mu_{\star} < a$ such that
$\underset{y \in M}{\inf}|D_rf_{\mu}(0,y)|>1 \quad \forall \mu \in (\mu_{\star},a)$.

\item $\exists$ \quad $0 < \alpha_1 < \alpha $ such that
$\|D_rf_{\mu}(r,y)\|_A<1 \quad \forall \mu
\in [0,a)$, where $A=\{ x\in \mathbb{R}^m : \alpha_1 \leq d(x,M) \leq \alpha \}$.

\item $\exists \quad \chi:[0,a) \rightarrow \mathbb{R}$ continuous with $0 \leq \chi(\mu)
\leq \alpha_1$ and $K(\mu):=\{x \in \mathbb{R}^m: \chi(\mu) \leq d(x,M) = d(x,y) \leq \alpha
\}$ such that $F_{\mu}(K(\mu)) \subseteq K(\mu) \quad \forall \mu \in (\mu_{\star},a)$.
Furthermore $c(\mu):=\|D_rf_{\mu}\|_{K(\mu)} < 1 \quad \forall \mu \in (\mu_{\star},a)$.

\item $c_{\star}(\mu):= (\|D_r f_{\mu}\|_{K(\mu)})
(1 + \|D_r \hat{g}_{\mu}\|_{K(\mu)}) + \|D_y f_{\mu}\|_{K(\mu)} <1$ for each $\mu \in
(\mu_{\star},a)$, where $\|.\|_{K(\mu)}$ is defined to be the $\sup$ norm over $K(\mu)$.
Here $(\hat{f}_{\mu},\hat{g}_{\mu})$ denotes the inverse map $F^{-1}_{\mu}$.

\item $(\|D_rf_{\mu}\|_{K(\mu)} + \|D_y f_{\mu}\|_{K(\mu)})
(\|D_r \hat{g}_{\mu}\|_{K(\mu)} + \|D_y \hat{g}_{\mu}\|_{K(\mu)}) \leq 1$ for each $\mu \in
(\mu_{\star}, a)$.

\item {\small $\sigma(\mu):= \|D_r f_{\mu}\|_{K(\mu)} (2\|D_r\hat{g}_{\mu}\|_{K(\mu)} +
\|D_y\hat{g}_{\mu}\|_{K(\mu)}) + \|D_yf_{\mu}\|_{K(\mu)}\|D_r\hat{g}_{\mu}\|_{K(\mu)} <
1$}
for all $\mu \in (\mu_{\star},a)$.
\end{enumerate}

Then for each $\mu \in (\mu_{\star},a)$, $\exists$ codimension-$1$ submanifolds $M_+(\mu)$ and
$M_-(\mu)$ in $K(\mu)$ such that both $M_+(\mu)$ and $M_-(\mu)$  are $F_{\mu}$-invariant, locally
attracting  and $\mathcal{C}^1$ diffeomorphic to $M$. $M$ is locally repelling, and $r>0$
for all $x=(r,y) \in M_+$, and $r<0$ for all $x=(r,y) \in M_-$.
\end{thm}
\begin{proof}
Recall that $N(\alpha)=\big \{x=(r,y) \in \mathbb{R}^m: |r|=d(x,M) \leq \alpha \quad, \quad y =
\pi (x)\big \}$. Now $ F_{\mu}(r,y)= (f_{\mu}(r,y), g_{\mu}(r,y))$ in component form
where, $f_{\mu}:N(\alpha) \rightarrow \mathbb{R}$ is the signed distance
 between $F_{\mu}(r,y)$ and $g_{\mu}(r,y)$ and $g_{\mu}:N(\alpha) \rightarrow M$ is the projection
 $\pi \circ F_{\mu}(y)$ of $F_{\mu}(r,y)$ on $M$. We shall break the proof up into a number
 of steps (claims).
\begin{claim}\label{cl:Mlocatt}$M$ is locally attracting  for $\mu \in (-a,0)$.
\end{claim}
\noindent
\emph{Proof of Claim \ref{cl:Mlocatt}:} Consider a point $(r_0, y_0) \in N(\alpha)$. Let
$(r_n,y_n)$
be the point obtained by applying the $n$-fold composition of $F_\mu$ with itself to
$(r_0,y_0)$. Then
\begin{displaymath}
(r_n,y_n)= F_{\mu}(r_{n-1},y_{n-1}) =
(f_{\mu}(r_{n-1},y_{n-1}),g_{\mu}(r_{n-1},y_{n-1}))
\end{displaymath}
implies that
\begin{displaymath}
d((r_n,y_n),M) = d((r_n,y_n),\pi(r_n,y_n))= f_{\mu}(r_{n-1},y_{n-1}).
\end{displaymath}
As $M$ is $F_{\mu}$-invariant, it follows that $f_{\mu}(0,y_{n-1})=0$ for all $n \in
\mathbb{N}$. \quad  So,
\begin{displaymath}
|r_n|=|f_{\mu}(r_{n-1},y_{n-1})|=|f_{\mu}(r_{n-1},y_{n-1}) - f_{\mu}(0,y_{n-1})|
= |\frac{\partial f_{\mu}(r^{\star},y_{n-1})}{\partial r}||r_{n-1}|
\end{displaymath}
by the mean value theorem. Thus $|r_n| <  c|r_{n-1}| < c^n|r_0|$, where\\
$c=\underset{(r,y) \in N(\alpha)}{\sup}|\frac{\partial f_{\mu}(r,y)}{\partial
r}|<1$ by property $(ii)$. Therefore, $r_n \rightarrow 0$ as $n \rightarrow \infty$.
Consequently $d((r_n,y_n),M) \rightarrow 0$.
That is, for any initial point $(r_0,y_0)$ in the neighborhood $N(\alpha)$ of $M$,
$F_{\mu}^n(r_0,y_0)$ converges to $M$. It follows that $M$ is locally attracting  for all
$\mu \in (-a,0)$.
\begin{claim}\label{cl:Mlocrep}$M$ is locally repelling for $\mu \in (\mu_{\star}, a)$.\end{claim}
\noindent
\emph{Proof of Claim \ref{cl:Mlocrep}:} Following the same steps as above, we find that
$|r_n| > c|r_{n-1}| > |r_{n-1}|$ whenever $|r_n|$ is sufficiently small owing to statement
$(iii)$. Accordingly the iterates
$\{x_n\}$ must eventually leave any sufficiently thin tubular neighborhood of $M$ for
$\mu \in (\mu_{\star}, a)$, which means that $M$ is locally repelling.

We now fix a $\mu \in (\mu_{\star},a)$ and suppress $\mu$ in the notation for simplicity. To
begin with, we shall prove the existence of $M_+$ as an $F_{\mu}$-invariant manifold
homeomorphic to $M$. It suffices to prove the existence of $M_+$, as the existence of $M_-$
can be established in the same way. Observe that $M_+$ is invariant iff $F(M_+) = M_+$.
We shall seek $M_+$ in the form of the graph of a continuous function over $M$ defined
as \begin{displaymath}M_+=\Gamma_{\psi} = \{(\psi(y),y):y \in M \},\end{displaymath}
where $M_+ \subset K = K(\mu)$, $\psi:M \rightarrow \mathbb{R}$ and $\psi(y) \geq 0$ for all $y \in
M$.
Then for all $y \in M$, we have that $(\psi(y),y) \in M_+$ iff $F(\psi(y),y)=(\psi(z),z) \in
M_+$, which is equivalent to
\begin{equation}\label{eq:frel}
\bigg(f \big(\psi(y),y \big),g\big(\psi(y),y\big)\bigg) =(\psi(z),z) \in M_+.
\end{equation}
$F$ is a diffeomorphism, hence $F^{-1}(\psi(z),z) = (\psi(y),y)$ which implies that
\begin{equation}\label{eq:finvrel}
\bigg(\hat{f}\big(\psi(z),z \big),\hat{g}\big(\psi(z),z\big)\bigg) = (\psi(y),y).
\end{equation}
where $F^{-1}=(\hat{f},\hat{g})$. Combining equations (\ref{eq:frel}) and
(\ref{eq:finvrel}), we find that $M_+$ is invariant iff $\psi$ satisfies the functional
equation
\begin{equation}\label{eq:funcrel}
\psi(z) = f\bigg(\psi\big(\hat{g}(\psi(z),z)\big), \hat{g}\big(\psi(z),z\big)\bigg).
\end{equation}
Let $Lip(A,B)$ denote the set of all Lipschitz functions from $A$ to $B$.
Let $\mathcal{L}(\psi)$ denote the Lipschitz constant of a Lipschitz function $\psi$,
and $\Gamma_{\psi}=\{(\psi(y),y): y \in M\}$ denote the
graph of $\psi$. Now define the set
\begin{displaymath}
X:=\{ \psi \in Lip(M,\mathbb{R}^+\cup \{0\} ): \mathcal{L}(\psi)\leq 1, \Gamma_{\psi}
\subseteq K\}.
\end{displaymath}
\begin{claim}\label{cl:XBanach}$\{X,\| \cdot \|_{K} \}$ is a Banach space.\end{claim}
\noindent
\emph{Proof of Claim \ref{cl:XBanach}:} Let $\{ \psi_n \}$ be a Cauchy sequence in $X$. Then for all $n \in
\mathbb{N}$ we have that $\psi_n : M \overset{Lip}{\longrightarrow} \mathbb{R}^+ \cup
\{0\}$, $\mathcal{L}(\psi_n) \leq 1$ and $\Gamma_{\psi_n} \in K(\mu)$. Here
$\Gamma_{\psi_n}$ denotes the graph of $\psi_n$ over $M$.  Since the sequence $\{ \psi_n \}$
is Lipschitz, every $\psi_n$ is continuous. Now $M$ is compact and $\mathbb{R}^+ \cup \{0\}$ is
closed, so the set of all continuous functions from $M$ to $\mathbb{R}^+\cup\{0\}$ with sup
norm forms a Banach space. Moreover, if $\psi_n \rightarrow \psi$ as $n \rightarrow \infty$,
it is clear that $\psi$ is also Lipschitz, with Lipschitz constant not greater than one.
Since $K$ is closed, $K$ contains all its limit points. Therefore
$\Gamma_{\psi_n} \in K$ for all $n$ implies that
\begin{displaymath}\underset{n\rightarrow \infty}{lim} \Gamma_{\psi_n}= \underset{n\rightarrow
\infty}{lim}(\psi_n(y),y) = (\psi(y),y) = \Gamma_{\psi} \in K,
\end{displaymath}
hence $X$ is a Banach space.
In view of  (\ref{eq:funcrel}), we define an operator  $\mathcal{F}$ on $X$ as follows.
\begin{equation}
\mathcal{F}(\psi)(y):=f\bigg(\psi\big(\hat{g}(\psi(y),y)\big),\hat{g}\big(\psi(y),y\big)\bigg).
\end{equation}
\begin{claim} \label{cl:Fcontraction}$\mathcal{F}(X) \subseteq X$.\end{claim}
\noindent
\emph{Proof of Claim \ref{cl:Fcontraction}:} Let $z=\hat{g}(\psi(y),y)$. Then $f(\psi(z),z)=$ signed distance
between $F(\psi(z),z)$ and $\pi(F(\psi(z),z))$. If $\psi(z) >0$, then $f(\psi(z),z)
>0$ since $F$ is side-preserving. So $\mathcal{F}(\psi)$ is indeed a function from $M$ to
$\mathbb{R}^+ \cup \{0\}$.
$\mathcal{F}(\psi)$ is continuous since it is a composition of continuous functions. Now it
follows from the mean value theorem and the definition of $X$ that
\begin{eqnarray*}
|\mathcal{F}(\psi)(y_1) - \mathcal{F}(\psi)(y_1)|&=& |f(\psi(z_1),z_1) - f(\psi(z_2),z_2)|\\
&\leq& \| \frac{\partial f}{\partial r} \|_{K} |\psi(z_1) - \psi(z_2)| + \|D_y f \|_{K}|z_1 - z_2| \\
&\leq& \bigg(\| \frac{\partial f}{\partial r} \|_{K} + \|D_y f\|_{K}\bigg) \quad |z_1 - z_2|
\end{eqnarray*}
where $\|.\|_K = sup\{|.|:(r,y) \in K \}$. The above inequality follows because $\psi$
is a Lipschitz function with Lipschitz constant $ \leq 1$. Also
\begin{eqnarray*}
|z_1 - z_2| &=& | \hat{g}(\psi(y_1),y_1) - \hat{g}(\psi(y_2),y_2)|\\
&\leq& \| D_r \hat{g}\|_K \quad |\psi(y_1) - \psi(y_2)| + \|D_y \hat{g} \|_K \quad |y_1 - y_2|\\
&\leq& \bigg( \| D_r \hat{g}\|_K + \|D_y g\|_K \bigg) \quad |y_1 - y_2|.
\end{eqnarray*}
The two inequalities obtained above, together with property $(vii)$ imply that
\begin{eqnarray*}
|\mathcal{F}(\psi)(y_1) - \mathcal{F}(\psi)(y_2)| &\leq& \bigg(\|\frac{\partial f}{\partial
r}\|_K + \| D_y f\|_K\bigg)\bigg(\|D_r \hat{g}\|_K + \|D_y \hat{g} \|_K\bigg)\\
&& \times |y_1 - y_1|\\
&\leq& |y_1 - y_2|.
\end{eqnarray*}
Therefore, $\mathcal{F}(\psi) \in Lip(M,\mathbb{R}^+\cup\{0\})$ and $\mathcal{L}(\mathcal{F}(\psi))
\leq 1$.
We will now prove that $\mathcal{F}$ is a contraction mapping.
Using statement $(vi)$ and the mean value theorem, we compute that
\begin{eqnarray*}
&&| \mathcal{F}(\psi_1)(y) - \mathcal{F}(\psi_2)(y) |\\
&&=|f(\psi_1(\hat{g}(\psi_1(y),y)),\hat{g}(\psi_1(y),y)) - f(\psi_2(\hat{g}(\psi_2(y),y)),\hat{g}(\psi_2(y),y)) |\\
&&\leq \left\|\frac{\partial f}{\partial r}\right\|_K |\psi_1(\hat{g}(\psi_1(y),y)) -
\psi_2(\hat{g}(\psi_2(y),y))|\\
&&\quad+ \|D_y f\|_K |\hat{g}(\psi_1(y),y) - \hat{g}(\psi_2(y),y)|\\
&&\leq  \left\|\frac{\partial f}{\partial r}\right\|_K |\psi_1(\hat{g}(\psi_1(y),y)) -
\psi_1(\hat{g}(\psi_2(y),y))|\\
&&\quad + \left\|\frac{\partial f}{\partial r}\right\|_K |\psi_1(\hat{g}(\psi_2(y),y)) -
\psi_2(\hat{g}(\psi_2(y),y))| \\
&& \quad + \|D_y f\|_K |\hat{g}(\psi_1(y),y) - \hat{g}(\psi_2(y),y)|\\
&& \leq \left\|\frac{\partial f}{\partial r}\right\|_K (\hat{g}(\psi_1(y),y) - \hat{g}(\psi_2(y),y) +
\|\psi_1 - \psi_2 \| )\\
&&\quad + \|D_y f\|_K  |\hat{g}(\psi_1(y),y) - \hat{g}(\psi_2(y),y)|\\
&&\leq \left\|\frac{\partial f}{\partial r}\right\|_K \big( \| D_r \hat{g}\|_K \| \psi_1 - \psi_2 \| +
\|\psi_1 - \psi_2 \|\big)+ \|D_y f\|_K
\|D_r \hat{g} \|_K \| \psi_1 - \psi_2 \|\\
&&= \left \{ \left\|\frac{\partial f}{\partial r}\right\|_K \big(1+\| D_r \hat{g} \|_K\big)
+ \|D_y f\|_K \right \} \|\psi_1 - \psi_2 \|.
\end{eqnarray*}
This is true for all $x \in K$. Hence it is true for the supremum with $x$ taken
over $K$; therefore, we obtain the relation
\begin{eqnarray*}
\| \mathcal{F}(\psi_1) - \mathcal{F}(\psi_2) \| &\leq& \left[ \|\frac{\partial f}{\partial
r}\|_K (1+\|D_r \hat{g} \|_K)+ \| D_yf \|_K \right]  \|\psi_1 - \psi_2 \| \\
\| \mathcal{F}(\psi_1) - \mathcal{F}(\psi_2) \| &\leq& c_{\star} \|\psi_1 - \psi_2 \|
\end{eqnarray*}
where $c_{\star}=\|\frac{\partial f}{\partial r}\|_K (1+\| D_r \hat{g} \|_K)+ \|D_y f\|_K$
is such that $0 < c_{\star} < 1$ by hypothesis $(vi)$. This implies that $\mathcal{F}(\psi) \in
Lip(M, \mathbb{R}^+\cup \{0\})$ and that $\mathcal{L}(\mathcal{F}(\psi)) < 1$. We note here
that the invariance of $M$ implies that $f(0,y)=\hat{f}(0,y)=0$ for all $(0,y) \in M$.
Accordingly $D_yf(0,y)=D_y\hat{f}(0,y)=0$ whenever $x=(0,y) \in M$, which means that both
$\|D_yf\|_K$ and $\|D_y\hat{f}\|_K$ can be made as small as we like by choosing a
sufficiently thin tubular neighborhood of $M$. Note that if $F(\psi(y),y)=(\psi(z),z)$,
then $(\psi(y),y)=F^{-1}(\psi(z),z)$, and it follows that
\begin{displaymath} \psi(y) = \hat{f}(\psi(z),z) \qquad {\rm and} \qquad y = \hat{g}(\psi(z),z).
\end{displaymath}
Now consider $\psi \in X$. By definition, we have
\begin{eqnarray*}
\Gamma_{\mathcal{F}(\psi)} &=& \{ (\mathcal{F}(\psi)(z),z) : z \in M\}\\
&=& \left \{ \Bigg(f\bigg(\psi(\hat{g}(\psi(z),z)), \hat{g}(\psi(z),z)\bigg),z\Bigg): z \in
M \right \}.
\end{eqnarray*}
This implies that
\begin{eqnarray*}
\Gamma_{\mathcal{F}(\psi)} &=& \{ \big(f(\psi(y),y),z\big): z \in M\}\\
&=& \{ \big(f(\psi(y),y),g(\psi(y),y)\big): g(\psi(y),y) \in M\}\\
&=& \{ \big(f(\psi(y),y),g(\psi(y),y)\big): y \in M\}\\
&=& \{ F(\psi(y),y) : y \in M \}.
\end{eqnarray*}
We know that $(\psi(y),y) \in K$ for all $y \in M$ and $F(K) \subseteq K$. This implies that
$\Gamma_{\mathcal{F}(\psi)} \subseteq K$, thereby proving the claim that $\mathcal{F}(X)
\subseteq X$. Hence $\mathcal{F}:X \rightarrow X$ is a contraction mapping with respect to
the sup norm on $X$.

Since $\mathcal{F}$ is a contraction on a complete metric space $X$, it has a unique fixed
point in $X$ owing to Banach's fixed point theorem. Let $\phi$ be the fixed point of
$\mathcal{F}$. Then $\phi \in Lip(M,\mathbb{R}^+\cup \{0\})$ with Lipschitz constant
$\mathcal{L}(\phi) \leq 1$, and $\phi$ satisfies the functional equation (\ref{eq:funcrel}).
Therefore,
\begin{equation}
\phi(z)=f\bigg(\phi\big(\hat{g}(\phi(z),z)\big),\hat{g}\big(\phi(z),z\big)\bigg).
\end{equation}
\begin{claim}\label{cl:M+exists} $M_+$ exists and is locally attracting. \end{claim}
\noindent
\emph{Proof of Claim \ref{cl:M+exists}:} We now define $M_+$ as the graph of $\phi$ as follows
\begin{displaymath} M_+= \Gamma_{\phi}= \{ (\phi(y),y); y \in M \},\end{displaymath}
where $\phi$ is as above. This proves the existence of $M_+$. That $M_+$ is locally attracting
follows directly from its definition as the graph of a fixed point (function) of a contraction
mapping. 
\begin{claim}\label{cl:M+homeo} $M_+$ is homeomorphic to $M$.\end{claim}
\noindent
\emph{Proof of Claim:} Let $H:M \rightarrow M_+$ be defined as $H(y):=(\phi(y),y)$. Then $H$
is injective, surjective and continuous. $H^{-1}$ exists and is also injective and
surjective (bijective). Since $M_+$ is compact, $H^{-1}$ is also continuous. 
Hence the manifold $M_+$ is homeomorphic to $M$.
\begin{claim}\label{cl:phiC1} The function $\phi$ is a class $\mathcal{C}^1$ map.\end{claim}
\noindent
\emph{Proof of Claim:} We know that $\phi$ is the solution to the functional equation
(\ref{eq:funcrel}), hence
$\phi(z)=f\bigg(\phi\big(\hat{g}(\phi(z),z)\big),\hat{g}\big(\phi(z),z\big)\bigg)$,
$\phi \in Lip(M,\mathbb{R}^+ \cup \{0\})$, and $\mathcal{L}(\phi) \leq 1$. We
will inductively construct a sequence of $\mathcal{C}^1$ functions $\psi_n$ which converges to
$\phi$. Then using the Arzela-Ascoli theorem, we will prove
that $\phi$ is $\mathcal{C}^1$. The details are as follows.

Choose $\psi_1$ to be a positive constant such that $\Gamma_{\psi} \subset K$. By construction,
$\psi_1$ is $\mathcal{C}^1$ and $\mathcal{L}(\psi_1) =0$. Now suppose $\psi_n$ is defined
and that $\psi_n$ is $\mathcal{C}^1$ with $\mathcal{L}(\psi_n) \leq 1$. We define
$\psi_{n+1}$ inductively as,
\begin{eqnarray}
\psi_{n+1}(z)&=& \mathcal{F}(\psi_n)(z)\\
\nonumber
&=& f\bigg(\psi_n\big(\hat{g}(\psi_n(z),z)\big),\hat{g}\big(\psi_n(z),z\big)\bigg).
\end{eqnarray}
Let $h_n:M \rightarrow \mathbb{R}^m$ denote the function $(\psi_n,Id)$, where $Id$ denotes
the identity map on the second coordinate. That is, $h_n(z)=(\psi_n(z),z)$. Then $h_n$ is
$\mathcal{C}^1$ by the induction hypothesis and the fact that it is the composition of
$\mathcal{C}^1$ maps, and
\begin{displaymath} \psi_{n+1}(z)=f \circ h_n \circ \hat{g} \circ h_n(z). \end{displaymath}
Here we have used the fact that both $f$ and $\hat{g}$ are $\mathcal{C}^1$ since $F$ and
$F^{-1}$ are $\mathcal{C}^1$ diffeomorphisms.

The sequence of functions $\{\psi_{n}(z)\}$ converges uniformly to $\phi(z)$ since $\phi$
satisfies the contractive functional equation (\ref{eq:funcrel}). The Jacobian of $\psi_{n+1}$
evaluated at $z$ is the $1 \times (m-1)$ matrix or the gradient vector of  $\psi_{n+1}$
 given as
\begin{equation}
D \psi_{n+1}(z)=Df\bigg(h_n\big(\hat{g}(h_n(z))\big)\bigg)Dh_n(\hat{g}(h_n(z)))D\hat{g}(h_n(z))Dh_n(z),
\end{equation}
owing to the chain rule. Moreover,
\begin{displaymath}\psi_{n+1}(z) = \mathcal{F}(\psi_{n})(z)\end{displaymath}
by construction. Hence, $\mathcal{L}(\psi_{n+1}) \leq 1$. Since, $\psi_{n+1}$ is
differentiable, this implies that $\|D\psi_{n+1}(z) \| \leq 1$ for all $z \in M$. By
induction, $\{D \psi _n (z)\}$ is a sequence of continuous functions, uniformly bounded by
$1$.

We will now prove the equicontinuity of $\{D \psi _n (z)\}$. The techniques used below are
actually global versions of the methods employed by Hartman \cite{Hartman} for local invariant manifolds,
and the role of the Lipschitz property follows an approach used by Hirsch
et al. \cite{HiPuSh} and Shub \cite{Shub} to study hyperbolic invariant manifolds.

For any function $\beta$, we define $\triangle \beta(z):= \beta(z + \triangle z) - \beta(z)$.
When $\triangle z$ is such that
$\triangle z \leq \min\{\delta, \frac{\delta}{\|D_r\hat{g}\|_K + \| D_y\hat{g}\|_K} \}$, we
will show that $ \| \triangle D\psi_n(z) \| \leq \tau(\delta)$
for all $n$, where $\tau$ depends only on $\delta$ and is such that $\tau(\delta)
\rightarrow 0$ as $\delta \rightarrow 0$. The desired result will be proved by induction as
follows.

For any $\delta >0$, we define quantities $\eta (\delta)$ and $\tau (\delta)$ as
\begin{eqnarray*}
&&\eta(\delta)=\sup\bigg\{\|D_rf(r + \triangle r, y + \triangle y) - D_rf(r,y)\|,\\
&& \|D_yf(r+ \triangle r, y + \triangle y) - D_yf(r,y)\|,
\|D_r\hat{f}(r + \triangle r, y + \triangle y) - D_r\hat{f}(r,y)\|,\\
&& \|D_y\hat{f}(r+ \triangle r, y + \triangle y) - D_y\hat{f}(r,y)\|,
\|D_rg(r+ \triangle r,y + \triangle y) - D_rg(r,y)\|,\\
&&\|D_yg(r+ \triangle r, y + \triangle y) - D_yg(r,y)\|,
\|D_r\hat{g}(r+ \triangle r, y + \triangle y) - D_r\hat{g}(r,y)\|,\\
&&\|D_y\hat{g}(r+ \triangle r, y + \triangle y) - D_y\hat{g}(r,y)\|
 :(r,y) \in N(\alpha), |\triangle r|, |\triangle y| \leq \delta \bigg\}
\end{eqnarray*}
and
\begin{eqnarray*}
\tau(\delta) &=& \frac{2(\|D_rf\|_K + \| D_y f \|_K + \|D_r\hat{g}\|_K + \|D_y\hat{g}\|_K)
}{1 - \sigma}\eta(\delta)
\end{eqnarray*}
where $\sigma <1$ is as defined in property $(viii)$. It is observed that $\eta (\delta)$
converges to $0$ as $\delta$ approaches $0$.

Recalling that $\psi_1$ is defined to be a constant, we have  $\triangle D\psi_1(z) \equiv
0$, and this implies that $\|\triangle D \psi_1(z) \| \leq \tau(\delta)$ for all $\triangle
z$. Suppose that $\|\triangle D\psi_n(z)\| \leq \tau(\delta)$ is satisfied whenever
$\triangle z
\leq$ min $\{\delta, \frac{\delta}{\|D_r\hat{g}\|_K + \| D_y\hat{g}\|_K} \}$.
Now,
\begin{displaymath} D\psi_{n+1}(z)= D \times C \times B \times A \end{displaymath}
where
\begin{displaymath}
D =\left[\begin{array}{cc} D_rf(\psi_n(\hat{g}(\psi_n(z),z)),\hat{g}(\psi_n(z),z))  &
D_yf(\psi_n(\hat{g}(\psi_n(z),z)),\hat{g}(\psi_n(z),z)) \end{array} \right]
\end{displaymath}
is a $1 \times m$ matrix,
\begin{displaymath}
C = \left[\begin{array}{c} D\psi_n(\hat{g}(\psi_n(z),z))\\ I_{m-1}\end{array} \right]_{m
\times (m-1)},
\end{displaymath}
\begin{displaymath} B = \left[\begin{array}{c} D\hat{g}(\psi_n(z),z)\end{array} \right]_{(m-1) \times m},
\end{displaymath}
and
\begin{displaymath} A = \left[\begin{array}{c}D\psi_n(z) \\ I_{m-1}\end{array} \right]_{m \times (m-1)}.
\end{displaymath}
Multiplying the four matrices above and taking into account the block matrix notation, it
follows that $D\psi_{n+1}(z)$ can be expressed in the following simpler form.
\begin{eqnarray*}
&&D\psi_{n+1}(z)\\
&&=D_rf(\psi_n(\hat{g}(\psi_n(z),z)),\hat{g}(\psi_n(z),z))
D\psi_n(\hat{g}(\psi_n(z),z))D_r\hat{g}(\psi_n(z),z)D\psi_n(z)\\
 &&\quad +
D_yf(\psi_n(\hat{g}(\psi_n(z),z)),\hat{g}(\psi_n(z),z))D_r\hat{g}(\psi_n(z),z)D\psi_n(z)\\
&& \quad + D_rf(\psi_n(\hat{g}(\psi_n(z),z)),\hat{g}(\psi_n(z),z))D\psi_n(\hat{g}(\psi_n(z),z))D_y\hat{g}(\psi_n(z),z)\\
&& \quad + D_yf(\psi_n(\hat{g}(\psi_n(z),z)),\hat{g}(\psi_n(z),z))D_y\hat{g}(\psi_n(z),z).
\end{eqnarray*}
Each of the four terms added above is a  $1 \times (m-1)$ vector. We will now estimate the
quantity $\| \triangle D \psi_{n+1}(z)\|$. Using the definitions, we find after a
straightforward calculation that $\triangle D\psi_{n+1}(z) = D\psi_{n+1}(z + \triangle z) - D
\psi_{n+1}(z)$ can be written in the form

\begin{eqnarray*}
&&\triangle D\psi_{n+1}(z)\\
&=&
\bigg\{D_r f(\psi_n(\hat{g}(\psi_n(z + \triangle z),z + \triangle z)),\hat{g}(\psi_n(z +
\triangle z),z + \triangle z))\\
&&\quad \times D\psi_n(\hat{g}(\psi_n(z + \triangle z),z + \triangle z))
D_r\hat{g}(\psi_n(z + \triangle z),z + \triangle z)\\
&&\quad \times D\psi_n(z + \triangle z)\\
&&- D_r f(\psi_n(\hat{g}(\psi_n(z),z)),\hat{g}(\psi_n(z),z))\\
&&\quad \times D\psi_n(\hat{g}(\psi_n(z),z))D_r\hat{g}(\psi_n(z),z)D\psi_n(z)\bigg\}\\
&&+ \bigg\{ D_yf(\psi_n(\hat{g}(\psi_n(z + \triangle z),z + \triangle z)),\hat{g}(\psi_n(z + \triangle z),z
+ \triangle z))\\
&&\quad \times D_r\hat{g}(\psi_n(z + \triangle z),z + \triangle z)D\psi_n(z + \triangle z)\\
&&- D_yf(\psi_n(\hat{g}(\psi_n(z),z)),\hat{g}(\psi_n(z),z))D_r\hat{g}(\psi_n(z),z)D\psi_n(z) \bigg\}\\
&&+ \bigg\{D_r f(\psi_n(\hat{g}(\psi_n(z + \triangle z),z + \triangle z)),\hat{g}(\psi_n(z + \triangle z),z +
\triangle z))\\
&&\quad \times D\psi_n(\hat{g}(\psi_n(z + \triangle z),z + \triangle z))D_y\hat{g}(\psi_n(z + \triangle z),z +
 \triangle z)\\
&&- D_r f(\psi_n(\hat{g}(\psi_n(z),z)),\hat{g}(\psi_n(z),z))D\psi_n(\hat{g}(\psi_n(z),z))D_y\hat{g}(\psi_n(z),z)\bigg\}\\
&&+ \bigg\{D_yf(\psi_n(\hat{g}(\psi_n(z + \triangle z),z + \triangle z)),\hat{g}(\psi_n(z + \triangle z),z
+ \triangle z))\\
&&\quad \times D_y\hat{g}(\psi_n(z + \triangle z),z + \triangle z)\\
&&- D_yf(\psi_n(\hat{g}(\psi_n(z),z)),\hat{g}(\psi_n(z),z))D_y\hat{g}(\psi_n(z),z)\bigg\}.
\end{eqnarray*}
We denote the four bracketed terms above as $T_1$, $T_2$, $T_3$ and $T_4$, respectively. For
instance,
\begin{eqnarray*}
T_2&=& D_yf(\psi_n(\hat{g}(\psi_n(z + \triangle z),z + \triangle z)),\hat{g}(\psi_n(z +
\triangle z),z + \triangle z))\\
&&\quad \quad \times D_r\hat{g}(\psi_n(z + \triangle z),z + \triangle z)D\psi_n(z + \triangle z)\\
&&\quad- D_yf(\psi_n(\hat{g}(\psi_n(z),z)),\hat{g}(\psi_n(z),z))D_r\hat{g}(\psi_n(z),z)D\psi_n(z).
\end{eqnarray*}
Adding and subtracting appropriate terms yields,
\begin{eqnarray*}
&&T_2\\
&&=\bigg\{ D_yf(\psi_n(\hat{g}(\psi_n(z + \triangle z),z + \triangle z)),\hat{g}(\psi_n(z +
\triangle z),z + \triangle z))\\
&&\quad \quad \times D_r\hat{g}(\psi_n(z + \triangle z),z + \triangle z)D\psi_n(z + \triangle z)\\
&&-D_yf(\psi_n(\hat{g}(\psi_n(z ),z)),\hat{g}(\psi_n(z),z ))D_r\hat{g}(\psi_n(z + \triangle
z),z + \triangle z)\\
&&\quad \quad \times D\psi_n(z + \triangle z)\bigg\}\\
&&+ \bigg\{D_yf(\psi_n(\hat{g}(\psi_n(z ),z)),\hat{g}(\psi_n(z),z ))D_r\hat{g}(\psi_n(z +
\triangle z),z + \triangle z)\\
&&\quad \quad \times D\psi_n(z + \triangle z)\\
&&\quad - D_yf(\psi_n(\hat{g}(\psi_n(z ),z)),\hat{g}(\psi_n(z),z ))D_r\hat{g}(\psi_n(z),z)
D\psi_n(z+\triangle z)\bigg\}\\
&& +
\bigg\{D_yf(\psi_n(\hat{g}(\psi_n(z ),z)),\hat{g}(\psi_n(z),z ))D_r\hat{g}(\psi_n(z),z)D\psi_n(z+\triangle z)\\
&& \quad - D_yf(\psi_n(\hat{g}(\psi_n(z),z)),\hat{g}(\psi_n(z),z))D_r\hat{g}(\psi_n(z),z)D\psi_n(z)\bigg\}.
\end{eqnarray*}
Using the triangle inequality and the above definition of the quantity $\eta(\delta)$, we
obtain
\begin{eqnarray*}
\|T_2\| &\leq& \eta(\delta)\|D_r\hat{g}\|_K \|D\psi_n\|_K + \|D_yf\|_K \eta(\delta)
\|D\psi_n\|_K\\
&& + \|D_yf\|_K\|D_r\hat{g}\|_K \|\triangle D\psi_n\|,
\end{eqnarray*}
where the first term is valid only when $\triangle z \leq \frac{\delta }{\|D_r\hat{g}\|_K +
\| D_y\hat{g}\|_K}$ and $\triangle z \leq \delta$. (The estimate on $\triangle z$ is
obtained by applying the chain rule and the mean value theorem to the first term in $T_2$.)
Since $\mathcal{L}(\psi_n) \leq 1$, it follows that $\|D\psi_n\| \leq 1$. By the induction
hypothesis, $\|\triangle D\psi_n\| \leq \tau(\delta)$. This implies that
\begin{equation}\label{eq:T2}
\|T_2\| \leq \eta(\delta)\bigg(\|D_r\hat{g}\|_K + \|D_yf\|_K\bigg) + \|D_yf\|_K\|D_r\hat{g}\|_K
\tau(\delta).
\end{equation}
Using similar analyses, we obtain
\begin{equation}\label{eq:T1}
\|T_1\| \leq \eta(\delta)(\|D_r\hat{g}\|_K + \|D_r f\|_K) + 2\|D_r f\|_K \|D_r\hat{g}\|_K
\tau(\delta),\\
\end{equation}
\begin{equation} \label{eq:T3}
\|T_3\| \leq \eta(\delta)(\|D_r f\|_K+ \|D_y\hat{g}\|_K ) +\|D_r f\|_K \|D_y\hat{g}\|_K
\tau(\delta)\\
\end{equation}
and
\begin{equation}\label{eq:T4}
\|T_4\| \leq \eta(\delta)(\|D_y\hat{g}\|_K + \|D_yf\|_K).
\end{equation}
Combining equations (\ref{eq:T2}) - (\ref{eq:T4}) gives,
\begin{eqnarray*}
\|\triangle D\psi_{n+1}\| &\leq 2 \left\{ \|D_r\hat{g}\|_K + \|D_y\hat{g}\|_K +
\|\frac{\partial
f}{\partial r}\|_K + \|D_yf\|_K \right \}\eta(\delta) \\
&+ \left \{\|\frac{\partial f}{\partial r}\|_K \big(2\|D_r\hat{g}\|_K + \|D_y\hat{g}\|_K
\big)+ \|D_yf\|_K\|D_r\hat{g}\|_K \tau(\delta) \right \}.
\end{eqnarray*}
Substituting the definition of $\sigma$ in the above inequality, we find that
\begin{displaymath} \|\triangle D\psi_{n+1}\| \leq \tau(\delta)(1-\sigma) + \sigma \tau(\delta),\end{displaymath}
which proves that
$ \|\triangle D\psi_{n+1}\| \leq \tau(\delta)$
whenever $\triangle z \leq$ min$\{ \delta, \frac{\delta}{\|D_r\hat{g}\|_K + \| D_y \hat{g}
\|_K} \}$.

Thus, we have proved by induction that $\|\triangle D\psi_n(z)\| \leq
\tau(\delta)$ (whenever $\triangle z$ is sufficiently small) for all $n$. The quantity
$\tau(\delta)$ is such that, $\tau(\delta) \rightarrow 0$ uniformly as $\delta$ approaches
$0$. Hence, the sequence of functions $\{D\psi_n(z)\}$ is equicontinuous.
Since the sequence $\{D \psi _n (z)\}$ is a uniformly bounded and equicontinuous
sequence of functions on a compact set $M$, it follows from the Arzela-Ascoli theorem
that there exists a subsequence $D\psi_{n_k}(z)$ which is
uniformly convergent on $M$. Let $\rho(z)$ be the uniform limit of $D\psi_{n_k}(z)$ as $k
\rightarrow \infty$. Since we know that $\psi_n$ converges to $\phi$, this implies that
$\rho = D \phi$. That is, $\phi$ is differentiable. Also, since $D \psi_n(z)$ is continuous
for every $n$ and the convergence is uniform, we find that $\rho$ is also continuous. That
is, $D \phi(z)$ is continuous. This implies that $\phi$ is class $\mathcal{C}^1$.

Hence, the map $H:M \rightarrow M_+$ defined earlier as $H(y):=(\phi(y),y)$ is a $\mathcal{C}^1$
diffeomorphism. Thus we have proved that the manifold $M_+$ is diffeomorphic to $M$.

Analogously, one can prove that $M_-$ is diffeomorphic to $M$. This proves that $M$ has
undergone a pitchfork bifurcation at $\mu_{\star}$, into a pair of locally attracting
invariant manifolds $M_+$ and $M_-$, each diffeomorphic to $M$, for each $\mu \in
(\mu_{\star},a)$. Thus the proof is complete.
\end{proof}

There is also a side-reversing version of Theorem \ref{thm:PBthethm} that can be proved in a
completely analogous manner; namely
\begin{thm}\label{thm:PBsiderev}
Let $F_{\mu}$ and $M$ satisfy all the hypotheses of Theorem \ref{thm:PBthethm}, except with
$F_{\mu}$ being side-reversing. Then for each $\mu \in (\mu_{\star},a)$, there exist
manifolds $M_-(\mu)$ and $M_+(\mu)$, both $\mathcal{C}^1$ diffeomorphic to $M$, such that
$F_{\mu}(M_+) = M_-$, $F_{\mu}(M_-) = M_+$ and $M_-(\mu) \cup M_+(\mu)$ is
$F_{\mu}$-invariant and locally attracting .
\end{thm}

In certain cases, the estimates in properties $(vi)$-$(viii)$ of Theorems \ref{thm:PBthethm}
and \ref{thm:PBsiderev} can be combined into a single statement, as in the following result.

\begin{cor}\label{cor:combinedthms} Let the hypotheses of Theorems \ref{thm:PBthethm} and
\ref{thm:PBsiderev} be as above, except that properties $(vi)$-$(viii)$ are replaced by the
single estimate\\
$(ix)$ \quad $\|D_rf_{\mu}\|_K \|D_r\hat{g}_{\mu}\|_K + (\|D_rf_{\mu}\|_K +
\|D_yf_{\mu}\|_K)\bigg(1+\|D_r \hat{g}_{\mu}\|_K\bigg)<1$ for each $\mu \in
(\mu_{\star},a)$. Then the conclusions of the theorem still follow.
\end{cor}
\begin{proof}One need only observe that $(vi)$-$(viii)$ follows directly from $(ix)$.
\end{proof}

\begin{rem} If the function $F_{\mu}$ and invariant manifold $M$ are of class
$\mathcal{C}^2$, property $(viii)$ of the above theorem is not essential, for in this case
the equicontinuity of the sequence $\{D\psi_n\}$ in the above proof follows from the mean
value theorem. Of course, if one wishes to prove the existence of bifurcated $\mathcal{C}^2$
diffeomorphs of $M$, an analog of property $(viii)$ involving second derivatives would be
necessary. Such an estimate, although rather complicated, can be obtained in a
straightforward manner, and we leave this to the reader. If both the map and invariant
submanifold are $\mathcal{C}^k$, with $k>2$, it is not difficult to obtain a $k$th order
derivative analog of $(viii)$ that would guarantee the existence of bifurcated
$\mathcal{C}^k$ diffeomorphs of $M$.
\end{rem}

\begin{cor}
Let the hypotheses be the same as in Theorems \ref{thm:PBthethm} and \ref{thm:PBsiderev}
with the following additional modifications: Property $(v)$ is replaced by\\
$(v')$ \quad $\exists \chi:[0,a)$ such that $0<\chi(\mu) \leq \alpha_1$, and
$F_{\mu}(K(\mu)) \subset K(\mu)$,\\ where $K(\mu)$ is as in $(v)$ for every $\mu_{\star} <\mu
< a$, and the following assumption is added.\\
$(x)$ \quad For every $\mu \in (\mu_{\star},a)$, $f_{\mu}(r,y) > r (<r)$ for $(r,y) \in
(0,\chi(\mu)]\times M$ and $F_{\mu}(r,y) <r (>r)$ for $(r,y) \in [-\chi(\mu),0) \times M$ in
the side-preserving (side-reversing) case.\\
Then, in addition to the conclusions of Theorem \ref{thm:PBthethm} and \ref{thm:PBsiderev},
we have the following dynamical properties: The submanifold $M_+(\mu)$ attracts all points
$x=(r,y) \in (0,\alpha]\times M$, and $M_-(\mu)$ attracts all points $x=(r,y)\in[\alpha,0)
\times M$ in the side-preserving case; and in the side-reversing case, $N(\alpha)\backslash
M$ is contained in the basin of attraction of $M_+(\mu) \cup M_-(\mu)$.
\end{cor}
\begin{proof}
We shall verify only the additional result for $M_+(\mu)$ in the side-reversing case, since
the proofs of all of the other cases are similar and require only obvious modifications. For
the case at hand, it obviously suffices to show that the iterates of a point $(r_0,y_0)$
with $0<r_0<\chi(\mu)$ eventually wind up in $K(\mu)$. Setting
$(r_n,y_n)=F_{\mu}^n(r_0,y_0)$, it follows from $(x)$ that $\{r_n\}$ is an increasing
sequence of real numbers, which must exceed $\chi(\mu)$ for $n$ sufficiently large. Thus the
proof is complete.
\end{proof}
\begin{rem}
It is natural to ask about the bifurcation phenomena that may occur when $0<\mu<\mu_{\star}$
and $M$ has regions where $|D_rf_{\mu}|>1$ and regions where $|D_rf_{\mu}|<1$. If $F_{\mu}$
leaves all points of $M$ fixed, one can readily prove the existence of \textquotedblleft
blistered\textquotedblright diffeomorphs of $M$ using one-dimensional theory. The
\textquotedblleft blister\textquotedblright regions, where $|D_rf_{\mu}|>1$, have a pair of
locally attracting copies of $M$ manifested as inner or outer blisters on $M$, while the
portion of $M$ inside the blister is locally repelling. However, when $F_{\mu}$ merely
leaves $M$ invariant without fixing all the points, the situation is much more complicated
and needs further investigation.
\end{rem}
\begin{rem}
A particularly useful feature of our main results, Theorem \ref{thm:PBthethm}, Theorem
\ref{thm:PBsiderev} ( and Theorem \ref{thm:PBctsthm} appearing in the sequel), is that they
are constructive. The desired bifurcated manifolds can be determined to any desired accuracy
by successive approximation. For example, to approximate $M_+(\mu)$ in the side-preserving
case, one simply starts with $\psi_1$ as a positive constant so that its graph is in
$K(\mu)$, and then computes successive approximations using the functional equation
(\ref{eq:funcrel}). The iterate $\psi_n$ for $n$ sufficiently large yields an approximation
$M_n$ that can be chosen to be arbitrarily $\mathcal{C}^1$ close to $M_+(\mu)$, and the
error can be estimated from the definition of the iterates.
\end{rem}
\section{Illustration of the (discrete) pitchfork bifurcation theorem}\label{sec:illus}
In this section, we illustrate Theorem \ref{thm:PBthethm} proved in Section \ref{sec:PBdiscrete}
with a canonical example.
Let $A \in SO_n(\mathbb{R})$, the special orthogonal group of real $n \times n$ matrices, comprised of
orthogonal matrices with determinant $1$. Define a linear map $L_A:\mathbb{R}^n \rightarrow
\mathbb{R}^n$ as
\begin{displaymath} L_A(x)=Ax.\end{displaymath}
The map $L_A$ is an analytic (linear) diffeomorphism. Every $(n-1)$-sphere $S_{\alpha}$ of radius $\alpha >0$ is
$L_A$-invariant. That is, $L_A(S_{\alpha}) = S_{\alpha}$
where $S_{\alpha}=\{x \in \mathbb{R}^n : |x|=\alpha \}$ and $\alpha >0$. Note that $S_1$
denotes a sphere of radius $1$, in the space on which $L_A$ acts. If for instance, $A \in
SO_2(\mathbb{R})$ then $L_A:\mathbb{R}^2 \rightarrow \mathbb{R}^2$ and $S_1$ is the same as
$S^1 \subset \mathbb{R}^2$. If $A \in SO_3(\mathbb{R})$ then $L_A:\mathbb{R}^3 \rightarrow \mathbb{R}^3$
and $S_1$ is the same as $S^2 \subset \mathbb{R}^3$. The subscript denotes the radius of the
sphere and the dimension of the sphere is one less than the ambient space.

\begin{figure}[!hbt]
\begin{minipage}{\fighalfwidth}
\includegraphics[width=\fighalfwidth]{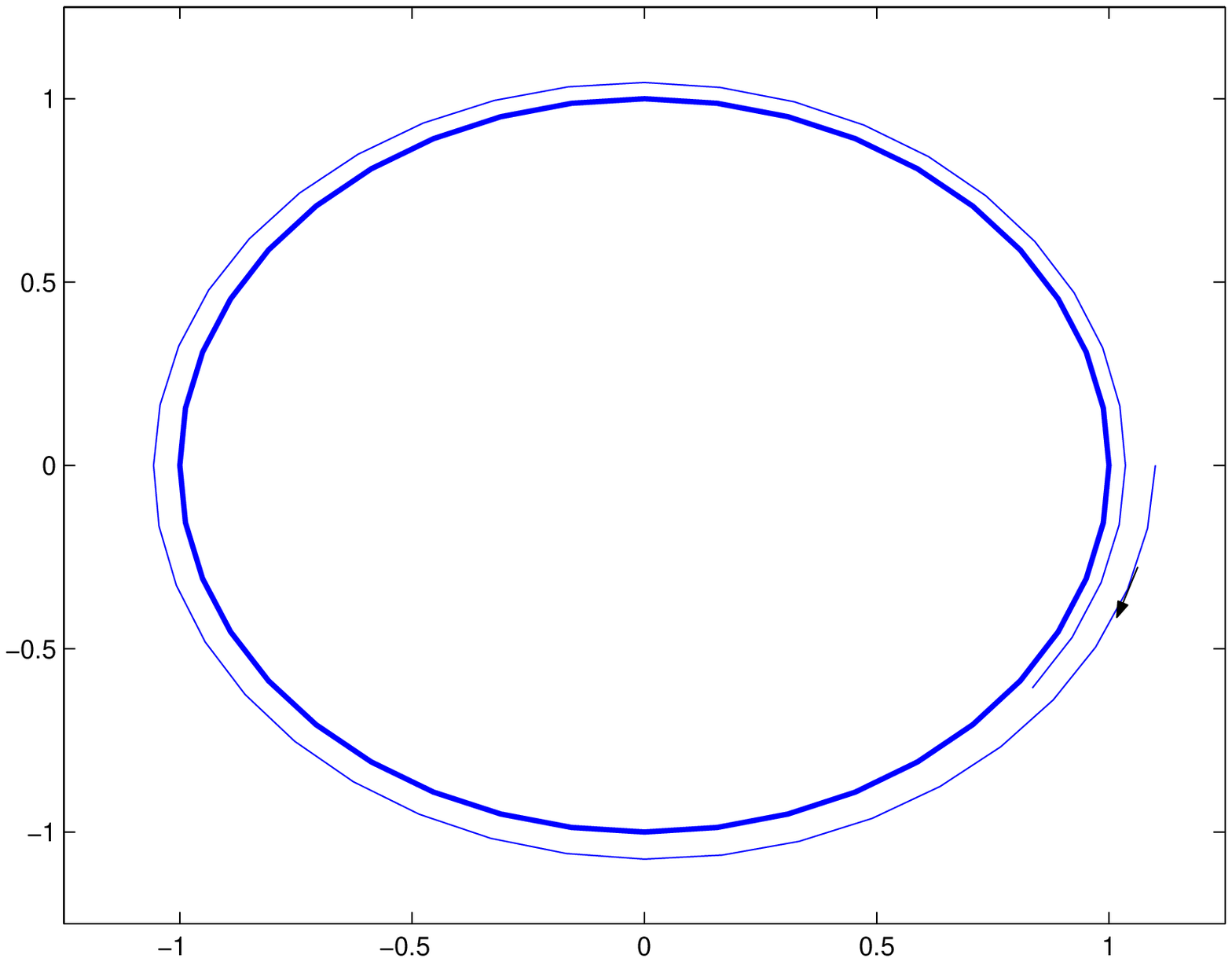}
\caption{Any trajectory outside $S_1$ converges to $S_1$.}
\label{fig:S1_02}
\end{minipage}
\quad
\begin{minipage}{\fighalfwidth}
\includegraphics[width=\fighalfwidth]{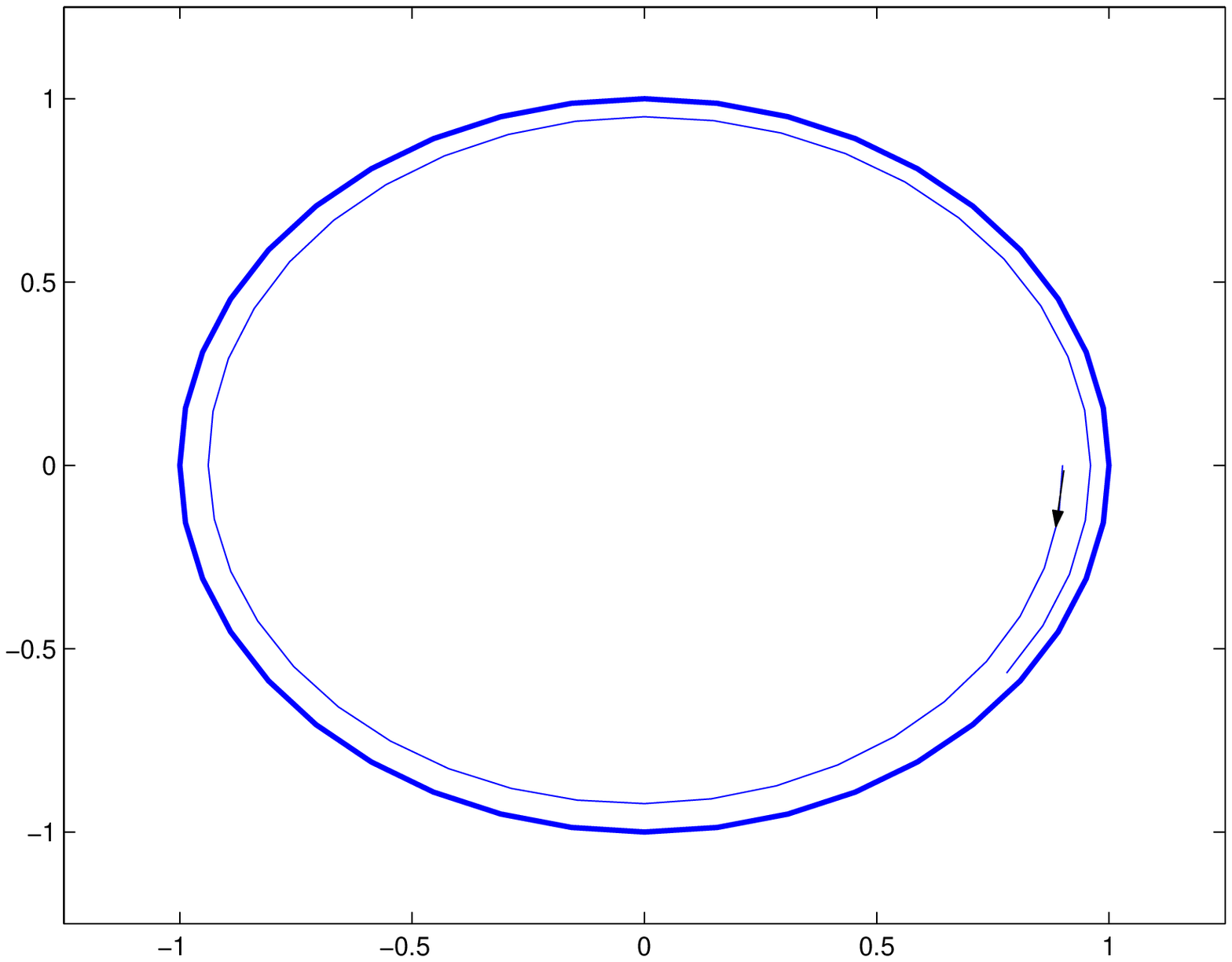}
\caption{Any trajectory inside $S_1$ converges to $S_1$.}
\label{fig:S1_03}
\end{minipage}
For $\mu<0$ that is, before the bifurcation $S_1$ is locally attracting.
Figures \ref{fig:S1_02} and \ref{fig:S1_03} are
generated for $\mu=\frac{-1}{50}$ in the canonical example.
\end{figure}

\begin{figure}[!hbt]
\begin{minipage}{\fighalfwidth}
\includegraphics[width=\fighalfwidth]{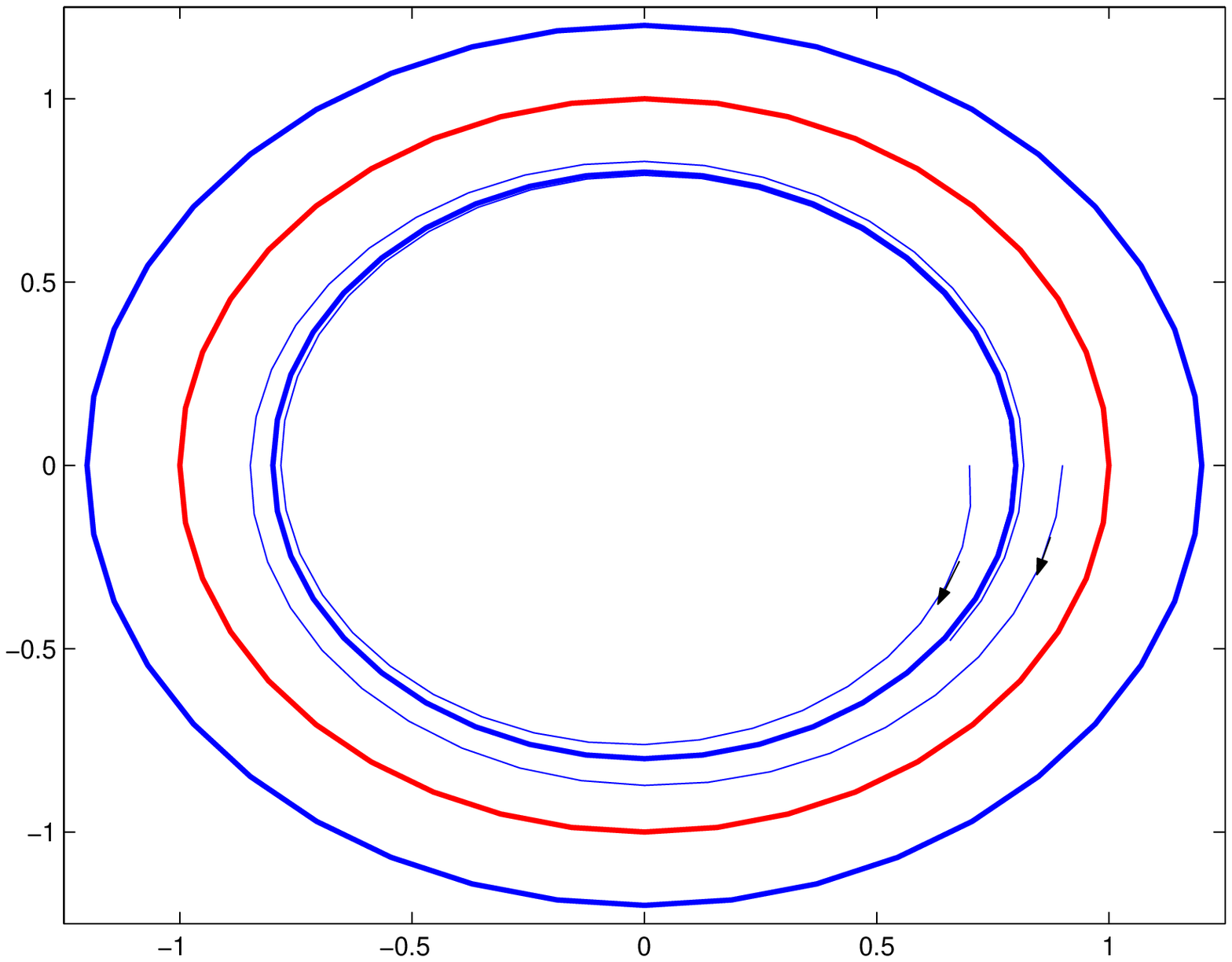}
\caption{Any trajectory outside converges to $S_{1+\sqrt{\mu}}$.}
\label{fig:S2_02}
\end{minipage}
\quad
\begin{minipage}{\fighalfwidth}
\includegraphics[width=\fighalfwidth]{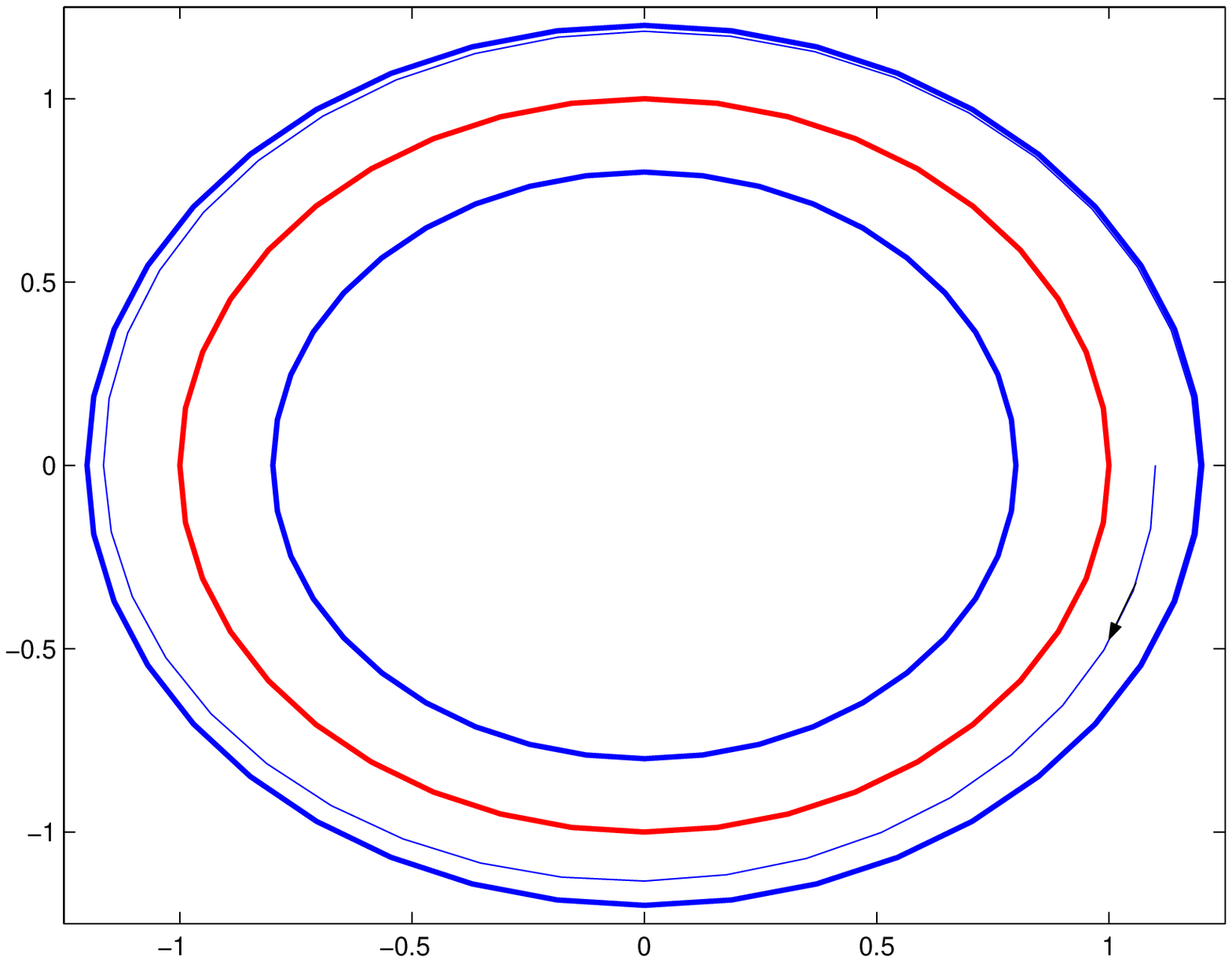}
\caption{Any trajectory inside converges to $S_{1-\sqrt{\mu}}$.}
\label{fig:S2_03}
\end{minipage}
For $\mu >0$, $S_1$ is locally repelling.
Figures \ref{fig:S2_02} and \ref{fig:S2_03} are generated for
$\mu=\frac{1}{50}$. Several points are iterated many times to generate these pictures.
\end{figure}

Now define $\sigma_{\mu}:[0,\infty) \rightarrow [0,\infty)$ to be a $\mathcal{C}^\infty$ function such
that $\sigma_{\mu}$ satisfies the following properties.
\begin{enumerate}
\item $\sigma_{\mu}' \equiv 0$ in a small neighborhood of $0$.
\item $\sigma_{\mu}(s) > 1$ for $0 \leq s < \frac{4}{5}$.
\item $\sigma_{\mu}(s)= 1- (s-1)^3 + \mu (s-1)$ for $\frac{4}{5} \leq s \leq \frac{6}{5}$.
\item $\sigma_{\mu}(s) <1$ for $\frac{6}{5} < s$.
\item $(s\sigma_{\mu}(s))' = s\sigma_{\mu}'(s) + \sigma_{\mu}(s) >0$ for $\mu \in
[\frac{-1}{25},\frac{1}{25}]$.
\end{enumerate}
We fix a matrix $A$ in $SO_{n}(\mathbb{R})$ and define $F_{\mu}:\mathbb{R}^n \rightarrow
\mathbb{R}^n$ as follows:
\begin{displaymath}F_{\mu}(x)=\sigma_{\mu}(|x|)L_A(x)= \sigma_{\mu}(|x|)Ax.\end{displaymath}
It is easy to see that $F_{\mu}$ is a diffeomorphism, and that it leaves $S_1$ invariant. That is,
$F_{\mu}(S_1)=S_1$.
The discrete dynamical system governed by $F_{\mu}$ is
\begin{equation}
x_{n+1}=F_{\mu}(x_n).
\end{equation}
In the notation of Section \ref{sec:PBdiscrete}, $M=S_1$. Due to the symmetry of the sphere
$S_1$, every point in $\mathbb{R}^n \backslash \{0\}$ can be uniquely described as being a radial
projection on $S_1$, so
the neighborhood $N(\alpha)$ is not restricted by the $\epsilon$-neighborhood theorem. However,
due to the nature of $\sigma_{\mu}$, we let $\alpha=\frac{1}{5}$ and consider the neighborhood
$N(\frac{1}{5})=\{x \in \mathbb{R}^n: |x|\in [\frac{4}{5},\frac{6}{5}]\}$.
We now check that all the hypotheses stated in Theorem \ref{thm:PBthethm} are satisfied.
\begin{enumerate}
\item Observe that  $F_{\mu}$ is side-preserving for $\mu \in [\frac{-1}{25},\frac{1}{25}]$ since $A$ preserves
orientation and $\sigma_{\mu}$ is positive-valued.\\
For this example,
\begin{displaymath}r=|x|-1 {\quad \rm and \quad}  y = \frac{x}{|x|}.\end{displaymath}

This implies that after a change of variables, $F_{\mu}(x)=\sigma_{\mu}(|x|)Ax$ becomes
\begin{eqnarray*}
F_{\mu}(r,y)&= \sigma_{\mu}(r+1)|x|A\frac{x}{|x|},\\
F_{\mu}(r,y)&= \sigma_{\mu}(r+1)(r+1)Ay.
\end{eqnarray*}
The property that $A$ preserves length is used in obtaining the above expression for
$F_{\mu}$, and again in finding $f_{\mu}$ and $g_{\mu}$ below:
\begin{eqnarray*}
f_{\mu}(r,y) = |\sigma_{\mu}(r+1) (r+1)Ay|-1 =(r+1)\sigma_{\mu}(r+1),\\
g_{\mu}(r,y)= Ay.
\end{eqnarray*}
This implies that
\begin{eqnarray*}
\frac{\partial f}{\partial r}=(r+1) \sigma_{\mu}'(r+1) + \sigma_{\mu}(r+1),\\
D_yf(r,y) \equiv 0, \quad D_rg(r,y) \equiv 0 \quad {\rm and} \quad D_yg(r,y) \equiv A.
\end{eqnarray*}
\item $\underset{(r,y)\in N(\frac{1}{5})}{\sup}|\frac{\partial f}{\partial r}|<1$ for all $\mu
\in [\frac{-1}{25},0)$ since the maximum $1$ is attained at $\mu=0$ as shown in Figure
\ref{fig:canonicalitem2}.

\begin{figure}[!hbt]
\begin{center}
\begin{minipage}{\fighalfwidth}
\includegraphics[width=\fighalfwidth,height=\fighalfheight]{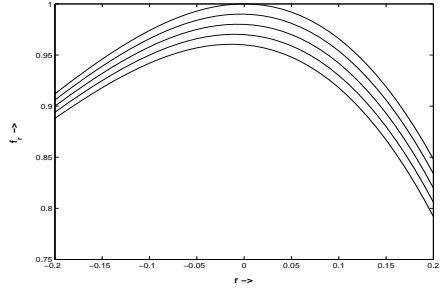}
\caption{$r$ vs $\frac{\partial f}{\partial r}$ for the canonical example for $r \in
[-0.2,0.2]$ as $\mu$ increases from $\frac{-1}{25}$ through $0$.}
\label{fig:canonicalitem2}
\end{minipage}
\end{center}
\end{figure}

\item For this example, $\mu_{\star}=0$ and $\inf |\frac{\partial f (0,y)}{\partial r}|>1$ for
all $\mu \in (0,\frac{1}{25}]$. The infimum is attained at $\mu=0$ as illustrated in Figure
\ref{fig:canonicalitem3}.

\begin{figure}[!hbt]
\begin{center}
\begin{minipage}{\fighalfwidth}
\begin{flushright}
\includegraphics[width=\fighalfwidth,height=\fighalfheight]{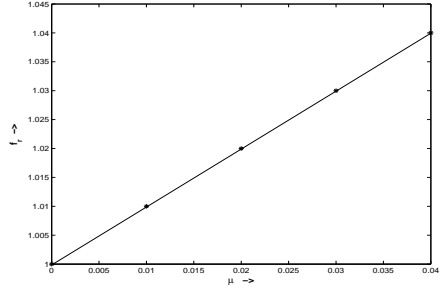}
\caption{Plot of $\mu$ vs $\frac{\partial f}{\partial r}$
for the canonical example for $r=0$ and $\mu$ in the interval $[0,\frac{1}{25}]$.}
\label{fig:canonicalitem3}
\end{flushright}
\end{minipage}
\end{center}
\end{figure}
\begin{figure}[!hbt]
\begin{center}
\begin{minipage}{\fighalfwidth}
\includegraphics[width=\fighalfwidth, height=\fighalfheight]{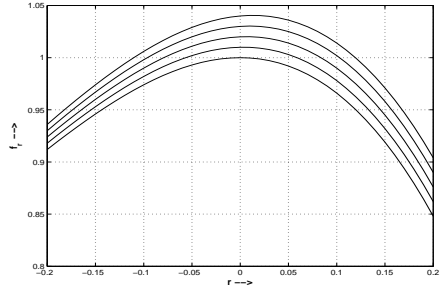}
\caption {Plot of $r$ vs $\frac{\partial f}{\partial r}$ for the canonical example  for $r \in [-0.2,0.2]$.}
\label{fig:canonicalitem4}
\end{minipage}
\end{center}
\end{figure}

\item For this case, $\alpha_1$ can be chosen to be $0.15$. As illustrated in Figure
\ref{fig:canonicalitem4},
$\underset{A}{\sup}|\frac{\partial f}{\partial r}|<1$, where $A=\{(r,y):0.15 \leq r \leq
0.2\}.$

\item $K(\mu)$ can be chosen to be $A$ and property $(v)$ follows from property $(iv)$.
\item Properties $(vi)$, $(vii)$ and $(viii)$ also follow from statement $(iv)$ since $D_rg_{\mu}(r,y) \equiv
0$, $D_yf(r,y) \equiv 0$ and $\|D_yg(r,y)\| =1$.
\end{enumerate}
Theorem  \ref{thm:PBthethm} implies that $S_1$ undergoes a pitchfork bifurcation at
$\mu_{\star}=0$. This is indeed the case and for $\mu \in (0,1/25]$: $F_{\mu}$ has three
invariant spheres $S_{1-\sqrt{\mu}}$, $S_1$ and $S_{1+\sqrt{\mu}}$ where $S_1$ is locally
repelling,
and $S_{1-\sqrt{\mu}}$ and $S_{1+\sqrt{\mu}}$ are locally attracting. This is illustrated in
Figures 
\ref{fig:S1_02}, \ref{fig:S1_03}, 
\ref{fig:S2_02} and \ref{fig:S2_03}.


\begin{rem}
The above example can be easily modified to illustrate Theorem \ref{thm:PBsiderev}.
Define a map $G_{\mu} = R\circ F_{\mu}$, where
$R:\mathbb{R}^m \backslash \{0\} \rightarrow \mathbb{R}^m \backslash \{0\}$ is a
smooth map such that $R(x)=\frac{2-|x|}{|x|}$ on the neighborhood $N(\frac{1}{5})$ of $S_1$.
Then $G_{\mu}$
is side-reversing, with all other properties the same as those of $F_{\mu}$. In this case,
$S_1$ undergoes a pitchfork bifurcation at $\mu_{\star} = 0$ and for $\mu \in
(0,\frac{1}{25}]$ the invariant manifolds are $S_1$ and $S_{1-\sqrt{\mu}} \cup
S_{1+\sqrt{\mu}}$. Note that $G_{\mu}(S_{1-\sqrt{\mu}}) = S_{1+\sqrt{\mu}}$ and $G_{\mu}(S_{1+\sqrt{\mu}})
= S_{1-\sqrt{\mu}}$.
\end{rem}
\section{Pitchfork bifurcation theorem for continuous dynamical system}
\label{sec:PBcts}
In this section, we state and prove a pitchfork bifurcation theorem for continuous dynamical
systems that is analogous to the result for the discrete case given by Theorem
\ref{thm:PBthethm}. The idea of the proof is to use the flow generated by a continuous
system to reduce the problem to the discrete system covered by Theorem \ref{thm:PBthethm}.
Consider a continuous dynamical system given by
\begin{equation}\label{eq:ctsdynsys}
\dot{x}=X(x,\mu)
\end{equation}
where $x\in \mathbb{R}^m$ and $\mu \in (-a,a) \subset \mathbb{R}$. The ($1$-parameter) vector field
$X(x,\mu)$, also denoted as $X_{\mu}(x)$, is assumed to be of class $\mathcal{C}^1$ in its domain.

Let $\phi(t,x,\mu)$, which we also denote by $\phi_{\mu}^t(x)$, be the unique solution (flow)
starting at $x$ when $t=0$, and $M$ be a compact, connected,
boundaryless, codimension-$1$, $\phi_{\mu}$-invariant manifold of $\mathbb{R}^m$, which means
$\phi_{\mu}^t(M)=M$ for all $|\mu|<a$. As in Section \ref{sec:PBdiscrete}, we define
\begin{displaymath} N(\alpha)=\{x\in \mathbb{R}^m: d(x,M) \leq \alpha\}\end{displaymath}
as a tubular neighborhood around $M$, where the $\epsilon$-neighborhood theorem \cite{GuPo}
is applicable. Any point $x$ in the region $N(\alpha)$ can be written as
\begin{displaymath} x=(r,y) \end{displaymath}
where $r$ is the signed distance between $x$ and the manifold $M$ and $y$ is the unique
point of $M$ closest to $x$. Recall that $r$ is positive if $x$ lies in the outer unbounded region of
$\mathbb{R}^m \backslash M$ and $r$ is negative if $x$ lies in the inner bounded region of
$\mathbb{R}^m \backslash M$. Again as in Section \ref{sec:PBdiscrete}, the notion of inner and outer
regions is obtained as an application of the Jordan-Brouwer separation theorem \cite{GuPo}.
 Note that $r=0$ when $x$ lies on $M$.

We shall assume that the vector field $X_{\mu}$ points into $N(\alpha)$ on $\partial
N(\alpha )$ for every $\mu$ in the set $(-a,a)$, which  means that positive semi-orbits of
(\ref{eq:ctsdynsys}) that begin in $N(\alpha)$ can never exit this tubular neighborhood. Note
that $X_{\mu}=(R_{\mu},Y_{\mu})$ in $r$,$y$-component form. Analogous to Section
\ref{sec:PBdiscrete}, it follows that $\phi_{\mu}^t$ maps $N(\alpha)$ into itself for all
$(t,\mu) \in [0,\infty) \times (-a,a)$. Now we can write the flow in terms of $r$ and $y$
components as
\begin{displaymath}\phi(t,(r,y),\mu)=(\rho(t,(r,y),\mu),\psi(t,(r,y),\mu)), \end{displaymath}
where $\phi(t,(r,y),\mu)$ is the signed distance between $\phi(t,(r,y),\mu)$ and $M$, and
$\psi(t,(r,y),\mu)$ is the normal projection of $\phi(t,(r,y),\mu)$ onto $M$.

In order to obtain the estimates necessary to reduce the continuous case to the discrete
case, we shall need to consider the derivative of the flow with respect to the initial
condition $x=(r,y)$. By $D\phi(t,x,\mu)$, we mean the Jacobian matrix of $\phi$ with respect
to $x$ defined as
\begin{displaymath}D\phi(t,x,\mu)=\left[
\begin{array}{cc}
D_r\rho(t,x,\mu) & D_y\rho(t,x,\mu)\\
D_r \psi(t,x,\mu) & D_y\psi(t,x,\mu)
\end{array}
\right].
\end{displaymath}
It is well known that the matrix $D\phi(t,x,\mu)$ is the unique solution of the initial
value problem (see e.g. in Hartman \cite{Hartman})
\begin{eqnarray}\label{eq:ctsIVP}
 \dot{\Phi}=D_xX_{\mu}(\phi(t,x,\mu))\Phi
 =\left[\begin{array}{cc}
D_rR_{\mu} & D_yR_{\mu}\\
D_r \psi_{\mu} & D_y\psi_{\mu}
\end{array}\right]_{\phi(t,x,\mu)}\Phi\\
\nonumber \Phi(0)=I_m,
\end{eqnarray}
where $I_m$ is the $m \times m$ identity matrix.

We shall make use of the following version of Gronwall's inequality.
\begin{lem}\label{lem:Gronwal}
Consider the linear matrix initial value problem
\begin{eqnarray}\label{eq:linGronwal}
\dot{\Phi}=\Gamma(t)\Phi\\
\nonumber \Phi(0)=I_m,
\end{eqnarray}
where $\Phi=(\phi_{ij})$, $\Gamma=(\gamma_{ij})$ and $\Gamma$ is a continuous matrix
function on the real line $\mathbb{R}$. Let $\phi(t)$, $\Phi_{I}(t)$, $\Phi_{II}(t)$,
$\Phi_{III}(t)$, $\gamma(t)$, $\Gamma_{I}(t)$, $\Gamma_{II}(t)$ and $\Gamma_{III}(t)$ be the
submatrices defined as follows:
\begin{displaymath} \phi(t)=\phi_{11}(t), \gamma(t)=\gamma_{11}(t),\end{displaymath}
\begin{displaymath}\Phi_{I}(t)=\left[ \phi_{12}(t), \cdots, \phi_{1m}(t)\right],
\Gamma_{I}(t)=\left[ \gamma_{12}(t), \cdots, \gamma_{1m}(t)\right],\end{displaymath}
\begin{displaymath}\Phi_{II}(t)=\left[ \phi_{21}(t), \cdots, \phi_{m1}(t)\right]^T,
\Gamma_{II}(t)=\left[ \gamma_{21}(t), \cdots, \gamma_{m1}(t)\right]^T,\end{displaymath}
\begin{displaymath}\Phi_{III}(t)=\left[
\begin{array}{ccc}\phi_{22}(t) & \cdots & \phi_{2m}(t)\\
\vdots & \cdots & \vdots\\
\phi_{m2}(t) & \cdots & \phi_{mm}(t)
\end{array}\right],
\Gamma_{III}(t)=\left[
\begin{array}{ccc}\gamma_{22}(t) & \cdots & \gamma_{2m}(t)\\
\vdots & \cdots & \vdots\\
\gamma_{m2}(t) & \cdots & \gamma_{mm}(t)
\end{array}\right],\end{displaymath}
where the superscript $T$ denotes transpose, so that
\begin{displaymath}\left[\begin{array}{cc}\dot{\phi} & \dot{\Phi}_{I}\\
\dot{\Phi}_{II} & \dot{\Phi}_{III}
\end{array}\right]=
\left[\begin{array}{cc}\gamma(t) & \Gamma_{I}(t)\\
\Gamma_{II}(t) & \Gamma_{III}(t)
\end{array}\right]
\left[ \begin{array}{cc}
\phi & \Phi_{I}\\
\Phi_{II} & \Phi_{III}
\end{array}\right]
\end{displaymath}
$\phi(0)=1$, $\Phi_{I}(0)=0$, $\Phi_{II}(0)=0$, $\Phi_{I}(0)=I_{m-1}$. Let $\sigma$, $\nu$ and $s$
be positive numbers satisfying
\begin{equation}\label{eq:ineq}
\sigma, \nu, \sigma^2, \nu^2 < s/4,
\end{equation}
and suppose that for some positive $t_{\star}$,
\begin{equation}\label{eq:gamestimates}
\gamma(t)\leq -2s, \quad  |\Gamma_{I}(t)|, |\Gamma_{II}(t)|\leq \sigma, \quad |\Gamma_{III}(t)|\leq \nu,
\end{equation}
whenever $|t| \leq t_{\star}$. Then for all $|t| \leq t_{\star}$ the solution of
(\ref{eq:linGronwal}) satisfies the estimates
\begin{eqnarray}\label{eq:phiestimates}
\nonumber
|\phi(t)| \leq E_0(t):= \kappa_1 e^{\lambda_-t} - \sigma(\lambda_+ + 2s)^{-1} \kappa_2
e^{\lambda_+ t},\\
\nonumber
|\Phi_{II}(t)| \leq E_2(t):= \sigma(\lambda -\nu)^{-1}\kappa_1 e^{\lambda_-t} - \kappa_2
e^{\lambda_+ t},\\
|\Phi_{I}(t)| \leq E_1(t):= \overline{\kappa}_1 e^{\lambda_-t} - \sigma(\lambda_+
2s)\overline{\kappa}_2e^{\lambda_+ t},\\
\nonumber
|\Phi_{III}(t)| \leq E_3(t):= \sigma(\lambda_- \nu)^{-1}\overline{\kappa}_1 e^{\lambda_-t} -
\overline{\kappa}_2e^{\lambda_+ t},
\end{eqnarray}
where
\begin{eqnarray}\label{eq:constants}
\lambda_{\pm}:=-\frac{(2s-\nu)}{2}\left[ 1 \pm \sqrt{1+ \frac{4\sigma^2+\nu^2}{(2s-\nu)^2}} \quad \right],\\
\nonumber
\kappa_1:=1-\sigma^2\{ (\lambda_+ +2s)[(\lambda_+ +2s) + \sigma^2(\lambda_- -\nu)]\}^{-1},\\
\nonumber
\kappa_2:=-\sigma\{ (\lambda_- \nu)(\lambda_+ +2s) [(\lambda_+ + 2s)+ \sigma^2(\lambda_- -\nu)]\}^{-1},\\
\nonumber
\overline{\kappa}_1:=\sigma(\lambda_- \nu)[(\lambda_+ +2s)(\lambda_- -\nu)+ \sigma^2]^{-1},\\
\nonumber
\overline{\kappa}_2:=\{1+\sigma^2[(\lambda_+ +2s)(\lambda_- -\nu)]^{-1}\}^{-1}.
\end{eqnarray}
\end{lem}
\begin{proof}
It follows from equations (\ref{eq:linGronwal})-(\ref{eq:gamestimates}) in the hypotheses that
\begin{displaymath}|\phi(t)| \leq u(t), \quad |\Phi_{I}(t)| \leq v(t), \quad |\Phi_{II}(t)| \leq w(t), \quad
|\Phi_{III}(t)| \leq z(t)
\end{displaymath}
for all $|t| \leq t_{\star}$, where $u$, $v$, $w$ and $z$ are the entries of the $2 \times
2$ matrix initial value problem
\begin{eqnarray}\label{eq:IVP}
\left[ \begin{array}{cc} \dot{u} & \dot{v}\\
\dot{w} & \dot{z} \end{array} \right] =
\left[ \begin{array}{cc} -2s & \sigma\\ \sigma & \nu \end{array} \right]
\left[ \begin{array}{cc} u & v\\ w & z \end{array} \right]\\
\nonumber
\left[ \begin{array}{cc} u(0) & v(0) \\ w(0) & z(0) \end{array} \right]=I_2
\end{eqnarray}
The eigenvalues - one negative and denoted by $\lambda_-$, and one positive and denoted by
$\lambda_+$ - of the constant matrix in (\ref{eq:IVP}) are easily computed and found to be
given by (\ref{eq:constants}). Now (\ref{eq:IVP}) can be solved by elementary means to yield
\begin{displaymath} u(t)=E_0(t), \quad v(t)=E_1(t), \quad w(t)=E_2(t), \quad z(t)=E_3(t), \end{displaymath}
where $E_0$, $E_1$, $E_2$ and $E_3$ are as defined in (\ref{eq:phiestimates}). Accordingly, we
have verified the desired estimates, thereby completing the proof.
\end{proof}
\begin{thm}\label{thm:PBctsthm} Let the vector field $X:\mathbb{R}^m \times(-a,a)$ be
$\mathcal{C}^1$, and let $M$ be a compact, connected, codimensions-$1$ invariant manifold
for (\ref{eq:ctsdynsys}) for every $\mu \in (-\sigma,\sigma)$. Suppose that the following properties
hold:
\begin{enumerate}
\item $X_{\mu}$ points into $N(\alpha)$ for all $(x,\mu)\in \partial N(\alpha) \times
(-a,a)$.
\item $D_rR(x,\mu) <0$ for all $x=(r,y)$ in the neighborhood $N(\alpha)$ for all $\mu \in
(-a,0)$.
\item There exists $0 \leq \mu_{\star} <a$ such that $D_rR(x,\mu) >0$ for all $(x,\mu) \in M
\times (\mu_{\star},a)$.
\item For each $\mu \in (\mu_{\star},a)$ there exists $0<\alpha_1(\mu) < \alpha$ and an
$s>0$ such that $X_{\mu}$ points into
\begin{displaymath} A(\mu)=\{ x \in \mathbb{R}^m: \alpha_1(\mu) \leq d(x,M) \leq \alpha \} \end{displaymath}
on its boundary and $D_rR((r,y),\mu) \leq -2s$ for $(r,y) \in A(\mu)$.
\item Let $\sigma$ and $\nu$ be positive constants such that $\sigma$, $\nu$, $\sigma^2$,
$\nu^2 < s/4$
\begin{displaymath}\|D_y R_{\mu}\|_{A(\mu)}, \quad \|D_rY_{\mu}\| \leq \sigma, \quad \|D_yY_{\mu}\|_{A(\mu)}
\leq \nu \end{displaymath}
for all $\mu \in (\mu_{\star},a)$, and $\sigma$, $\nu$ are sufficiently small with respect to
$s$ so that
\begin{equation}\label{eq:Eest1}
E_0(t)(1+E_2(-t)) + E_1(t)<1,
\end{equation}
\begin{equation}\label{eq:Eest2}
(E_0(t)+ E_1(t))(E_2(-t) + E_3(-t))\leq 1,
\end{equation}
\begin{equation}\label{eq:Eest3}
E_0(t)(2E_2(-t) +E_3(-t))+ E_1(t)E_2(-t))< 1,
\end{equation}
where $E_0$, $E_1$, $E_2$ and $E_3$ are as in Lemma \ref{lem:Gronwal}, for all $\mu \in
(\mu_{\star},a)$ and each $1 \leq t \leq 2$. Then the invariant submanifold $M$ is locally
attracting for $\mu \in (-a,0)$, and locally repelling for $\mu \in (\mu_{\star},a)$.
Furthermore, for each $\mu \in (\mu_{\star},a)$ there exist a pair of $\mathcal{C}^1$
diffeomorphs $M_+(\mu)$ and $M_-(\mu)$ of $M$ in $A(\mu)$ such that both $M_+(\mu)$ and
$M_-(\mu)$ are invariant for (\ref{eq:ctsdynsys}) and locally attracting.
\end{enumerate}
\end{thm}
\begin{proof}Using the relation $(\dot{r}, \dot{y})=(R((r,y),\mu), Y((r,y),\mu))$, it
follows from the mean value theorem that
\begin{eqnarray*}
\dot{r}&=R((r,y).\mu)=R((r,y).\mu) - R((0,y),\mu)\\
&= D_rR((r_{\star},y),\mu)r,
\end{eqnarray*}
where $r_{\star}$lies between $0$ and $r$, and we have used the property $R((0,y),\mu)=0$,
which follows from the invariance of $M$. Consequently, property $(ii)$ implies that
$\dot{r} < 0$ when $r>0$ and $\dot{r} > 0$ when $r<0$, which means that trajectories tend
toward $M$ as $t$ increases. Hence, $M$ is locally attracting for each $-a < \mu <0$.
Similarly, one can also use the mean value theorem to show that it follows from property
$(iii)$ that $M$ is locally repelling for $\mu \in (\mu_{\star},a)$.

From here on, we fix $\mu \in (\mu_{\star},a)$ and suppress it in order to simplify the
notation. We shall first show that for each $t \in [1,2]$, the map $T^t$ defined as
\begin{displaymath}T^t(x):=\phi(t,x),\end{displaymath}
where $\phi$ is the flow generated by the differential equation (\ref{eq:ctsdynsys}),
satisfies the hypotheses of Theorem \ref{thm:PBthethm}. As $\Phi=D\phi$ satisfies the
initial value problem (\ref{eq:ctsIVP}), Lemma \ref{lem:Gronwal} is an ideal instrument for
proving the desired result.

Observe that from the form of the estimates (\ref{eq:phiestimates}) of Lemma
\ref{lem:Gronwal}, that for a given $s>0$ it is indeed possible to select positive numbers
$\sigma$, $\nu$ sufficiently small for estimates (\ref{eq:Eest1})-(\ref{eq:Eest3}) of
property ($iv$) to obtain for all $1 \leq t \leq 2$. This is with the understanding that we
may assume without loss of generality that we are in $T^2(A)$, so that we can still take
advantage of the initial norm estimates for the terms with arguments $-t$ in
(\ref{eq:Eest1})-(\ref{eq:Eest3}), which correspond to the inverse of $T^t$ owing to the
group property of the flow $\phi$. Accordingly it follows from Lemma \ref{lem:Gronwal} and
Theorem \ref{thm:PBthethm} that for each $\mu \in (\mu_{\star},a)$ and $t \in [1,2]$, the
map $T^t$ has a unique pair of contractive invariant manifolds $M^t_{\pm}$ in $A$, which are
$\mathcal{C}^1$-diffeomorphic with $M$.

We shall now show that the manifolds $M^t_+$ and $M^t_-$ are, in fact, the same for all
$t\in \mathbb{R}$, and they are invariant for the entire flow $\phi$. It is enough to verify
this for $M^t_+$, as the proof for $M^t_-$ is identical. Consider any rational number of the
form $q=1+\frac{1}{m}$ lying between $1$ and $2$. Then, by definition
\begin{displaymath}T^q(M^q_+)=M^q_+.\end{displaymath}
Applying the map $T^q$, $n$ times to this equation yields
\begin{displaymath}[T^q]^n(M^q_+)=M^q_+,\end{displaymath}
which by the additivity property of the flow becomes
\begin{displaymath}T^{(n+m)}(M^q_+)=M^q_+.\end{displaymath}
But the unique contractive manifold for $T^1$ is $M^1_+$, and $T^{(n+m)}$ is an $(n+m)$-fold
composite of $T^1$ with itself. Hence,
\begin{displaymath}T^{(n+m)}(M^1_+)=M^1_+,\end{displaymath}
so it follows from uniqueness that
\begin{equation} \label{eq:M+uniq}
M^q_+=M^1_+,
\end{equation}
and this must hold for all rational numbers $1\leq q \leq 2$.

It now follows from (\ref{eq:M+uniq}), the completeness of real numbers, and the continuity of
the flow that $M^t_+=M^1_+$ for all $t\in [1,2]$, and
\begin{equation}\label{eq:M+uniqin12}
T^t(M^1_+)=M^1_+
\end{equation}
for all $1\leq t \leq 2$. Any $t>2$ can be written as $t=m+\tau$, where $m$ is a positive
integer and $\tau \in [1,2]$. Consequently,
\begin{eqnarray*}
T^t(M^1_+)&=T^{(m+\tau)}(M^1_+)=T^m\circ T^{\tau}(M^1_+)\\
&= T^m(M^1_+)=M^1_+,
\end{eqnarray*}
owing to (\ref{eq:M+uniqin12}) and the definition of $M^1_+$. Whence (\ref{eq:M+uniqin12})
holds for all $t\geq 1$. In fact, it holds for all $|t| \geq 1$ since
\begin{displaymath} T^{-t}(T^t(M^1_+))=M^1_+=T^{-t}(M^1_+) \end{displaymath}
whenever $t\geq 1$.

Finally, for any $\epsilon >0$,
\begin{displaymath}
T^{\epsilon}(M^1_+)=T^{(1+\epsilon)}\circ T^{-1}(M^1_+)=T^{(1+\epsilon)}(M^1_+)=M^1_+.
\end{displaymath}
Thus $T^t(M^1_+)=M^1_+$ for all $t \geq 0$, and it therefore follows as above that the same
is true for all $t<0$ as well. We denote the unique invariant manifold $M^1_+$ by $M_+$.
This yields the desired result that $T^t(M_+)=M_+$ for all $t\in \mathbb{R}$, which completes
the proof.
\end{proof}

\begin{rem} For continuous dynamical systems, only the side-preserving case can occur.
\end{rem}

Just as in the case of a discrete dynamical system, we can obtain a more complete
description of the dynamical systems in $N(\alpha)$ by making a minor additional assumption.
\begin{cor} In addition to the hypotheses of Theorem \ref{thm:PBctsthm}, suppose that
$R((r,y),\mu)$ is positive (negative) for $0<r\leq \alpha_1(\mu)$ ($-\alpha_1(\mu) \leq r
<0$) whenever $\mu \in (\mu_{\star},a)$. Then $M_+(\mu)$ attracts all points of $N(\alpha)$
with $r>0$ and $M_-(\mu)$ attracts all points of $N(\alpha)$ with $r<0$.
\end{cor}
\begin{proof}
The additional property guarantees that the positive semi-orbits in $N(\alpha)$ not lying in
$M$, eventually enter $A(\mu)$ and are then attracted to $M_+(\mu)$ or $M_-(\mu)$. This
completes the proof.
\end{proof}

\section{Conclusions}
We have proved that codimension-$1$, compact invariant manifolds in discrete dynamical systems,
 undergo pitchfork bifurcations when the system satisfies suitable conditions. The hypotheses of
 the theorem are
easily verifiable estimates on the norms of partial derivatives of the function determining the discrete
dynamical system, which makes this result well suited to a variety of applications. When the bifurcation
parameter $\mu$ is between $0$ and $\mu_{\star}$, some portions of $M$ may be locally
repelling and some locally attracting (in the normal direction), so the proof of our theorem
would need to be modified to handle this case, which is an interesting subject for future
investigation.

The case when the whole manifold $M$ bifurcates into $M_-$ and $M_+$ as $\mu$
increases through zero, corresponds to $\mu_{\star}=0$. The fact that $\mu_{\star}$ can
be greater than $0$ allows for $M$ to eventually bifurcate and does not impose the
restriction that $M$ bifurcate all at once. The theorem is slightly weaker than the theorem
in one-dimension since the theorem does not completely determine the dynamics of the system in the
region between a neighborhood of $M$ and the neighborhood $A$ of $M_-$ and $M_+$.

The pitchfork bifurcation in $\mathbb{R}$ is assumed to be one stable fixed point bifurcating
into two stable fixed points separated by an unstable fixed point. We have generalized this
result to a compact, connected, boundaryless, codimension-$1$, locally attracting invariant
submanifold of $\mathbb{R}^m$ that becomes locally repelling and bifurcates into two locally
attracting diffeomorphic copies of itself separated by the locally repelling manifold. The
techniques we have used here should enable us to obtain new results on higher dimensional
versions of other types of bifurcations such as Hopf and saddle-node Hopf bifurcations (see
e.g. \cite{KrOl}). We plan to investigate these and related generalizations in the future.

\end{document}